
\font\addressfont=cmr8
\newcount\s
\newcount\n
\def\clearn{\n =0}
\def\N{\global\advance\n by 1 \global\s=1
                       \global\clearnn
                       \global\clearnnn
                       \global\clearnnnn
                       {\the\n}%
                     }
\newcount\nn
\def\clearnn{\nn =0}
\def\NN{\global\advance\nn by 1 \global\s=2
                        \global\clearnnn
                        \global\clearnnnn
                        {\the\n}.{\the\nn}%
                   }
\newcount\nnn
\def\clearnnn{\nnn =0}
\def\NNN{\global\advance\nnn by 1 \global\s=3
                        \global\clearnnnn
                        {\the\n}.{\the\nn}.{\the\nnn}%
                     }

\newcount\nnnn
\def\clearnnnn{\nnnn =0}
\def\NNNN{\global\advance\nnnn by 1 \global\s=4
                        {\the\n}.{\the\nn}.{\the\nnn}.{\the\nnnn}%
                     }

\newif\ifindi
\newbox\boxindice
\inditrue
\newwrite\ind
\def\apriindice{\ifindi%
                            \setbox\boxindice=\vbox{\input \nomefile.ind }%
                            \immediate\openout\ind=\nomefile.ind%
                            \immediate\write\ind{{\sectionfont
Contents}\vskip10pt}%
                            \else\fi}
\def\indice#1{\ifindi\immediate\write\ind{{#1}}\else\fi}


\newif\ifrifer
\rifertrue
\newwrite\rifer
\def\apririferimenti{\ifrifer%
                            \immediate\openout\rifer=\nomefile.rif%
                            \immediate\write\rifer{riferimenti}%
                            \else\fi}
\def\riferimenti#1{\ifrifer\immediate\write\rifer{#1}\else\fi}
\def\chiudiriferimenti{\ifrifer\immediate\closeout\rifer\else\fi}


\def\rif#1{\expandafter\xdef\csname s#1\endcsname{\number\s}%
            \ifnum\s>0 \expandafter\xdef\csname n#1\endcsname{\number\n}\fi%
            \ifnum\s>1 \expandafter\xdef\csname nn#1\endcsname{\number\nn}\fi%
            \ifnum\s>2 \expandafter\xdef\csname nnn#1\endcsname{\number\nnn}\fi%
            \ifnum\s>3 \expandafter\xdef\csname nnnn#1\endcsname{\number\nnnn}\fi%
            \riferimenti{\cite{#1}........#1}%
                   }

\def\cite#1{%
            \expandafter\ifcase\csname s#1\endcsname%
            \or \csname n#1\endcsname%
            \or \csname n#1\endcsname.\csname nn#1\endcsname%
            \or \csname n#1\endcsname.\csname nn#1\endcsname.\csname nnn#1\endcsname%
            \or \csname n#1\endcsname.\csname nn#1\endcsname.\csname nnn#1\endcsname.%
                \csname nnnn#1\endcsname%
                   \fi}

\font \proc=cmbx10
\def\proclaim#1. #2\par
                  {\medbreak  {\proc #1. \enspace }{\sl #2\par }
                    \ifdim \lastskip <\medskipamount \removelastskip
                            \penalty 55\medskip \fi}

\def\biblio#1#2{ \item{\hbox to 1.2cm{[#1]\hss}}
                           {#2} }


\def\numfont{\bf}

\def\prop#1#2{\proclaim{{\numfont\csname#1\endcsname.} Proposition}. {#2}\par}
\def\theo#1#2{\proclaim{{\numfont\csname#1\endcsname.} Theorem}. {#2}\par}
\def\lemm#1#2{\proclaim{{\numfont\csname#1\endcsname.} Lemma}. {#2}\par}
\def\defi#1#2{\proclaim{{\numfont\csname#1\endcsname.} Definition}. {#2}\par}
\def\coro#1#2{\proclaim{{\numfont\csname#1\endcsname.} Corollary}. {#2}\par}

\long\def\proo#1#2{{{\numfont\csname#1\endcsname.}\proc \ Proof.}
                                    #2
                                    \hfill$\square$\vskip 10pt
                                }
\long\def\exam#1#2{{{\numfont\csname#1\endcsname.}\proc \ Example.}
                                    #2
                                    \vskip 10pt
                                }

\long\def\rema#1#2{{{\numfont\csname#1\endcsname.}\proc \ Remark.}
                                    #2
                                    \vskip 10pt
                                }

\long\def\nota#1#2{{{\numfont\csname#1\endcsname.}\proc \ Notations.}
                                    #2
                                    \vskip 10pt
                                }

\long\def\proof#1{{\proc \ Proof.}
                                    #1
                                    \hfill$\square$\vskip 10pt
                                }


\font\titlefont=cmbx10 scaled \magstep2
\font\authorfont=cmcsc10 scaled \magstep1
\font\sectionfont=cmbx10 scaled \magstep1
\font\subsectionfont=cmbx10
\font\abstractfont=cmr8
\font\titleabstractfont=cmcsc8

\def\section#1{\vskip 10pt
                          {\sectionfont\N.\ #1}
                          \indice{\proc \the\n.\ #1\endgraf}
                          \vskip 10pt}
\def\subsection#1{\vskip 8pt
                             {\subsectionfont\NN.\ #1}
                             \indice{\quad\proc \the\nn.\ #1\endgraf}
                             \vskip 8pt}
\def\intro{\vskip 10pt
                          {\sectionfont Introduction.}
                          \vskip 10pt\indice{\proc Introduction.\endgraf}}
\def\refer{\vskip 10pt
                          {\sectionfont References.}
                          \vskip 10pt\indice{\proc References.\endgraf}}
\def\title#1{\vskip 5pt
                          \centerline{\titlefont #1}
                          \vskip 10pt}
\def\author#1{\vskip 8pt
                          \centerline{\authorfont #1}
                          \vskip 10pt}

\def\abstract#1{
\centerline{\vtop{\hsize=12truecm\baselineskip=9pt\strut
                 {\titleabstractfont Abstract. }
                         \abstractfont #1}}
                          \vskip 10pt}

\clearn

\vsize= 230truemm
\hsize= 155truemm
\voffset= 8truemm
\hoffset= 3truemm

\abovedisplayskip=6pt plus2pt minus 4pt
\belowdisplayskip=6pt plus2pt minus 4pt
\abovedisplayshortskip=0pt plus2pt
\belowdisplayshortskip=2pt plus1pt minus 1pt

\input amssym.def
\input amssym.tex

\def\point{{\scriptscriptstyle\bullet}}
\def\lra{\longrightarrow}
\def\ra{\rightarrow}

\def\Hom{{\rm Hom}}
\def\tot{{\rm tot}}

\bigskip 


\input amssym.def
\input amssym.tex
\input epsf

\def\abstra{ \ifnum\pageno<0\else\pageno=-1\fi
          \vfill\eject
          \ifodd\pageno{}\else\quad\eject\fi
          \global\inichapter=\pageno
          \xdef\chaptertitle{Abstract}
          \xdef\sectiontitle{Abstract}
          \line{\hfill\chaptertitlefont Abstract}
          \vskip 30pt
           }


\def\oggi{\number\day\space\ifcase\month\or
       gennaio\or febbraio\or marzo\or aprile\or maggio\or giugno\or
       luglio\or agosto\or settembre\or ottobre\or novembre\or dicembre\fi
        \ \number\year}

\def\today{\ifcase\month\or
               January\or February\or March\or April\or May\or June\or
               July\or August\or September\or October\or November\or December\fi
               \space\number\day, \number\year}

\def\aujourdhui{\number\day\space\ifcase\month\or
               Janvier\or F\'evrier\or Mars\or Avril\or May\or Juin\or
               Juillet\or Ao\^ut\or Septembre\or Octobre\or Novembre\or 
               D\'ecembre\fi
               \ \number\year}


\catcode`\@=11   
\long\def\ulap#1{\vbox to \z@{\vss#1}}
\long\def\dlap#1{\vbox to \z@{#1\vss}}

\def\xlap#1{\hbox to \z@{\hss#1\hss}}
\long\def\ylap#1{\vbox to \z@{\vss#1\vss}}


\def\ulap#1{\vbox to \z@{\vss#1}}
\def\dlap#1{\vbox to \z@{#1\vss}}
\def\hlap#1{\hbox to \z@{\hss#1\hss}}
\def\vlap#1{\vbox to \z@{\vss\hbox{#1}\vss}}
\def\clap#1{\vbox to \z@{\vss\hbox to \z@{\hss#1\hss}\vss}}

\catcode`\@=12   

\catcode`\@=11

\def\mmatrix#1{\null\,\vcenter{\normalbaselines\m@th
                    \ialign{\hfil$##$\hfil&&\ \hfil$##$\hfil\crcr
                    \mathstrut\crcr\noalign{\kern-\baselineskip }
                    #1\crcr\mathstrut\crcr\noalign{\kern-\baselineskip}}}\,}

\catcode`\@=12

\catcode`\@=11

\def\Rightarrowfill{ $\m@th \mathord =\mkern -6mu\cleaders
       \hbox {$\mkern -2mu\mathord =\mkern -2mu$}
          \hfill \mkern -6mu\mathord \Rightarrow $}
\def\Leftarrowfill{$\m@th \mathord \Leftarrow \mkern -6mu\cleaders
             \hbox {$\mkern -2mu\mathord =\mkern -2mu$}
             \hfill \mkern -6mu\mathord =$}
\def\Leftrightarrowfill{$\m@th \mathord \Leftarrow \mkern -6mu\cleaders
             \hbox {$\mkern -2mu\mathord =\mkern -2mu$}
             \hfill \mkern -6mu\mathord \Rightarrow $}

\def\leftrightarrowfill{$\m@th \mathord \leftarrow \mkern -6mu\cleaders
             \hbox {$\mkern -2mu\mathord -\mkern -2mu$}
             \hfill \mkern -6mu\mathord \rightarrow $}

\def\underrightarrow#1{\vtop {\m@th \ialign {##\crcr
                     $\hfil \displaystyle { #1 }\hfil $\crcr
                     \noalign {\kern 1pt \nointerlineskip }
                    \rightarrowfill \crcr \noalign {\kern 0pt }}}}

\def\underleftarrow#1{\vtop {\m@th \ialign {##\crcr
                     $\hfil \displaystyle { #1 }\hfil $\crcr
                     \noalign {\kern 1pt \nointerlineskip }
                    \leftarrowfill \crcr \noalign {\kern 0pt }}}}

\def\underleftrightarrow#1{\vtop {\m@th \ialign {##\crcr
                     $\hfil \displaystyle { #1 }\hfil $\crcr
                     \noalign {\kern 0pt \nointerlineskip }
                    \leftrightarrowfill \crcr \noalign {\kern 0pt }}}}

\def\overleftrightarrow#1{\vbox {\m@th \ialign {##\crcr
                   \noalign {\kern 3pt } \leftrightarrowfill \crcr
                   \noalign {\kern 3pt \nointerlineskip }
                   $\hfil \displaystyle { #1 }\hfil$\crcr }}}

\catcode`\@=12

\def\rightto#1{\hbox to#1pt{\rightarrowfill}}
\def\leftto#1{\hbox to#1pt{\leftarrowfill}}
\def\leftrightto#1{\hbox to#1pt{\leftrightarrowfill}}

\def\mapstoto#1{\mapstochar\mkern -4mu\hbox to#1pt{\rightarrowfill}}
\def\rmapstoto#1{\hbox to#1pt{\leftarrowfill}\mkern -4mu\mapstochar}

\def\hookrightto#1{\lhook\mkern -8mu\hbox to#1pt{\rightarrowfill}}
\def\hookleftto#1{\hbox to#1pt{\leftarrowfill}\mkern -8mu\rhook}

\def\twoheadrightto#1{\hbox to#1pt{\rightarrowfill}\mkern -20mu\rightarrow }
\def\twoheadleftto#1{\leftarrow\mkern -20mu\hbox to#1pt{\leftarrowfill}}

\def\hooktwoheadrightto#1{\lhook\mkern -8mu%
                               \hbox to#1pt{\rightarrowfill}%
                               \mkern -20mu\rightarrow }
\def\hooktwoheadleftto#1{\leftarrow\mkern -20mu%
                              \hbox to#1pt{\leftarrowfill}%
                              \mkern -8mu\rhook}




\newdimen\rotdimen
\def\vspec#1{\special{ps:#1}}
\def\rotstart#1{\vspec{gsave currentpoint currentpoint translate
        #1 neg exch neg exch translate}}
\def\rotfinish{\vspec{currentpoint grestore moveto}}

\def\rotr#1{\rotdimen=\ht#1\advance\rotdimen by\dp#1%
        \hbox to\rotdimen{\hskip\ht#1%
        \vbox to\wd#1{\rotstart{90 rotate}%
        \box#1\vss}\hss}\rotfinish}%

\def\rotr#1{\rotdimen=\ht#1\advance\rotdimen by\dp#1%
        \hbox to\rotdimen{\hskip\ht#1
        $\vcenter to\wd#1{\rotstart{90 rotate}%
        \box#1\vss}$\hss}\rotfinish}%

\def\rotl#1{\rotdimen=\ht#1\advance\rotdimen by\dp#1%
        \hbox to\rotdimen{\vbox to\wd#1{\vskip\wd#1\rotstart{270 rotate}%
        \box#1\vss}\hss}\rotfinish}%

\newbox\horarr
\def\Uh#1{\setbox\horarr=\hbox{$\hookleftto{#1}$}\rotr{\horarr}}
\def\Dh#1{\setbox\horarr=\hbox{$\hookrightto{#1}$}\rotr{\horarr}}
\def\Ut#1{\setbox\horarr=\hbox{$\twoheadrightto{#1}$}\rotl{\horarr}}
\def\Dt#1{\setbox\horarr=\hbox{$\twoheadleftto{#1}$}\rotl{\horarr}}
\def\Uht#1{\setbox\horarr=\hbox{$\hooktwoheadleftto{#1}$}\rotr{\horarr}}
\def\Dht#1{\setbox\horarr=\hbox{$\hooktwoheadrightto{#1}$}\rotr{\horarr}}
\def\Uto#1{\setbox\horarr=\hbox{$\mapstpto{#1}$}\rotl{\horarr}}
\def\Dto#1{\setbox\horarr=\hbox{$\rmapstoto{#1}$}\rotl{\horarr}}



\def\cdef#1#2#3{
        \vbox to0pt {\vss \hbox to0pt {\hss #3 \hss} \vss}
}
\newbox\deform
\newdimen\altezza
\newdimen\profondita
\newdimen\larghezza
\def\cdef#1#2#3{
        \setbox\deform=\hbox{#3}
        \altezza=\ht\deform%
        \profondita=\dp\deform%
        \larghezza=\wd\deform%
       \ifnum #1>0 \multiply\larghezza by #1 \else\multiply\larghezza by -#1\fi%
       \ifnum #2>0 \multiply\altezza by #2 \else\multiply\altezza by -#2\fi%
       \ifnum #2>0 \multiply\profondita by #2 \else\multiply\profondita
by -#2\fi%
       \advance\altezza by \profondita%
        \hbox to\larghezza{\hss\lower\profondita
        \vbox to \altezza{\vss{
        \vbox to0pt {\vss \hbox to0pt {\hss #3 \hss} \vss}
        \vss}}
        \hss}
        }


\catcode`\@=11
\def\smalleqalign#1{\null \,\vcenter{\openup \jot \m@th %
       \ialign {\strut \hfil $\scriptstyle {##}$%
                &$\scriptstyle {{}##}$\hfil \crcr #1\crcr }}\,}
\def\smallmatrix#1{\null \,\vcenter {\baselineskip=6pt \m@th 
             \ialign {\hfil $\scriptstyle ##$\hfil &&\  
             \hfil $\scriptstyle ##$\hfil \crcr 
             \mathstrut \crcr \noalign {\kern -7pt } #1\crcr 
             \mathstrut \crcr \noalign {\kern -7pt }}}\,}
\def\verysmallmatrix#1{\null \,\vcenter {\baselineskip=4pt \m@th 
             \ialign {\hfil $\scriptscriptstyle ##$\hfil &&\  
             \hfil $\scriptscriptstyle ##$\hfil \crcr 
             \mathstrut \crcr \noalign {\kern -5pt } #1\crcr 
             \mathstrut \crcr \noalign {\kern -5pt }}}\,}

\catcode`\@=12


\def\ss{\scriptstyle}
\def\sss{\scriptscriptstyle}

\def\ov{\overline}

\def\phi{\varphi}

\def\theta{\vartheta}

\def\rho{\varrho}

\def\epsilon{\varepsilon}

\def\NoBlackBoxes{\global\overfullrule 0pt}
\NoBlackBoxes

\def\ra{\mathop{\rightarrow}\limits}

\def\lra{\mathop{\longrightarrow}\limits}


\def\leq{\leqslant}
\def\geq{\geqslant}




\font\tenscr=rsfs10 
\font\sevenscr=rsfs7 
\font\fivescr=rsfs5 
\skewchar\tenscr='177 \skewchar\sevenscr='177 \skewchar\fivescr='177
\newfam\scrfam \textfont\scrfam=\tenscr \scriptfont\scrfam=\sevenscr
\scriptscriptfont\scrfam=\fivescr
\def\scr{\fam\scrfam}

\def\b#1{ {\Bbb {#1}} }

\def\c#1{ {\cal #1} }
\def\c#1{ {\scr #1} }

\def\cOmega{{\mit\Omega}}
\def\cTheta{{\mit\Theta}}


\def\Hom{{\rm Hom}}

\def\cHom{{{\c H}{\mit o m}}}

\def\O{{\rm O}}

\def\bydef{\mathrel{:=}}


\def\limpr{\mathop{\underleftarrow{\rm lim}}}

\def\limin{\mathop{\underrightarrow{\rm lim}}}

\def\limpro{\mathop{\lim\limits_{\displaystyle\leftarrow}}}

\def\limind{\mathop{\lim\limits_{\displaystyle\rightarrow}}}



\def\Pro{{\rm Pro}}
\def\Ind{{\rm Ind}}


\def\bydef{\mathrel{:=}}

\def\Diff{{\rm Diff}}

\def\Strat{{\rm Strat}}

\def\epsilon{\varepsilon}

\def\point{{\scriptscriptstyle\bullet}}
\def\DR{{\rm DR}}
\def\tDR{\widetilde{\rm DR}}

\def\Q{{\rm Q}}

\def\bR{{\bf R}}

\def\tot{{\rm tot}}

\def\circdot{{\circ\kern-3.5pt\lower .8pt\hbox{$\cdot$}}}


\abovedisplayskip=6pt plus2pt minus 4pt
\belowdisplayskip=6pt plus2pt minus 4pt
\abovedisplayshortskip=0pt plus2pt
\belowdisplayshortskip=2pt plus1pt minus 1pt


\def\I(#1){{\rm Ind(\c #1)}}
\def\P(#1){{ \rm Pro(\c #1)}}
\def\bs{{\bigskip}}
\def\O(#1){{\rm {Ob}(\c #1)}}
\def\o(#1){{\rm {Ob}(\rm #1)}}
\def\Coh{\rm {Coh}}

\def\QCoh{\rm {QCoh}}

\def\opp{\circ}

\input xy
\xyoption{all}

\def\nomefile{ProStrat}
\title{Differential Complexes and Stratified Pro-Modules}
\author{Luisa Fiorot}

\abstract { 
In this paper we introduce the category of stratified Pro-modules 
and the notion of induced object in this category.
We propose a translation of Morihiko Saito equivalence results ([S.2])
using the dual language of Pro-objects. 
So we prove an equivalence between the derived category of 
stratified Pro-modules and the category of Pro-differential complexes.
We also supply a comparison with the notion of Crystal in Pro-module
(introduced by P. Deligne in 1960).
		}
		\bigskip
		{\addressfont
{
Luisa Fiorot

Universit\`a degli Studi di Padova 

Dipartimento di Matematica Pura ed Applicata	

Via Trieste, 63 - 35121 Padova

Italy 

fiorot@math.unipd.it}}
\bigskip


\apririferimenti

\bigskip
\intro{
Let $X$ be a smooth separated noetherian scheme of finite type 
over $\Bbb C$. 
We first note that a differential operator of finite order $m\in \Bbb N$
(between $\c O_X$-modules) can be defined
in two different ways: the first using induced right $\c D_{X,m}$-modules,
as done by Morihiko Saito in [S.2], the second using
 the sheaf of principal parts $\c P^m_{X}$.
Thus, any differential complex 
$\c L^\point$
admits two ``linearized" versions.
The first is given by M. Saito's functor 
$\DR^{-1}_X(\c L^\point )\bydef \c L^\point\otimes_{\c O_X}\c D_X$ 
in the category of right $\c D_X$-modules while
the other is given by Grothendieck's formalization functor 
$\Q^0_X\bydef \{\c P^m_X\}_{\b Z}\otimes_{\c O_X}\c L^\point$
in the category of stratified Pro-modules.
\endgraf
In [S.2] M. Saito proved the equivalence between
the derived category of right $\c D_X$-modules (quasi-coherent
as $\c O_X$-ones) and a suitable
localized category of differential complexes.
Our aim is to prove a ``dual" version of this equivalence
replacing quasi-coherent right $\c D_X$-modules by 
Pro-coherent stratified ones.
The main idea is that of using the Grothendieck formalization functor
${\Q}^0_{X}$ instead of Saito's $\DR^{-1}_{X}$.
On the other hand, a functor $\DR_X$ is always defined
on stratified objects simply by taking horizontal sections.
Suitably localizing these functors
gives an equivalence between the derived category
of Stratified Pro-coherent modules and that of
Pro-differential complexes (suitably localized).
We also define a category of $\cOmega^\point_X$-modules in
Pro-object (suitable localized) and we prove an equivalence 
with that of stratified Pro-modules
(as done in the dual case in [F1]).
\endgraf
In the last section we interpret stratified Pro-coherent modules 
as objects in the Crystalline site.
In particular we prove that
the category of stratified Pro-coherent modules
is equivalent to the category of ``crystals in Pro-modules".
This last notion was first introduced by P. Deligne in 
a cycle of lectures he
gave at IHES.
There Deligne proposed the notion of
``Crystals in Pro-modules"  attached to algebraic
constructible sheaves on an analytic space $X^{\ss an}$.
Moreover he proved that the category of ``regular crystals in Pro-modules"
is equivalent to that of ``algebraic" constructible sheaves
on $X^{\ss an}$ (unfortunately this work was not published).
By Deligne's equivalence theorem we obtain the
equivalence between the derived category of
regular stratified Pro-coherent modules 
and that of ``algebraic" constructible sheaves, and thus
a sort of Riemann-Hilbert correspondence.
\smallskip
As noted above the notion of
stratified Pro-coherent module is dual to
that of quasi-coherent right $\c D$-module.
In a work in progress
we expect to prove an anti-equivalence of categories between the 
category of perfect $\c D$-complexes and that of perfect complexes  of
stratified Pro-modules.
This anti-equivalence is compatible with the duality in the category
of differential complexes.
In particular when any object of a differential
complex is  coherent on $\c O_X$, the notion of 
$\c  D_X$-qis (see [S.2]) is equivalent to that of $Q^0_X$-qis.
\endgraf

\smallskip
I would like to thank Pierre Berthelot for his suggestions and for his notes on
Deligne lectures on ``Cristaux discontinues". 
I would also like to thank  Maurizio Cailotto and Anne Virrion for the improvements 
they suggested to me during the preparation of this work.
}\bs
\section{Pro-Coherent $\c O_X$-Modules}
We briefly recall some results on the category of Pro-coherent $\c O_X$-modules.
\defi{NN}
{([SGAIII,1]).
By definition the category 
$\Pro(\Coh(\c O_X))\bydef \Ind(\Coh(\c O_X)^{\opp})^\opp$.
Objects are filtering projective systems of coherent $\c O_X$-modules, while morphisms between
two such objects
$\c F_I$, $\c G_J$ are $\Hom_{\Pro}(\c F_I,\c G_J)=\limpro_J\limind_I\Hom_{\c O_X}(\c F_i,\c G_j)$.
For brevity we will use  the notation 
$\nu(\c O_X)$ for $\Pro(\Coh(\c O_X))$
and
$\mu(\c O_X)$ for $\Ind(\Coh(\c O_X))$.
} 
\smallskip

\rema{NN}
{\rif{rem1}
Given  a noetherian scheme $X$ over $\Bbb C$,
the category $\mu(\c O_X)$ is equivalent to the category of quasi-coherent
$\c O_X$-modules (denoted by $\QCoh(\c O_X)$) see [RD, appendix].
Moreover any object in $\mu(\c O_X)$
may be represented as an inductive system whose transition
morphisms  are injective maps.
}
\smallskip
\rema{NN}
{
The functor tensor product:
$$\matrix{\_\otimes_{\c O_X}\_:&
\Coh(\c O_X)\times \Coh(\c O_X)&\lra& \Coh(\c O_X)\cr
&\hfill (\c F,\c G) & \longmapsto &\c F\otimes_{\c O_X}\c G\cr}
$$
extends to the procategory:
$$\matrix{\_\otimes_{\c O_X}\_:&
\Pro(Coh(\c O_X))\times \Pro(\Coh(\c O_X))&\lra&\Pro(\Coh(\c O_X))\cr
&\hfill (\{\c F_i\}_I,\{\c G_j\}_J) & \longmapsto &
\{\c F_i\otimes_{\c O_X}\c G_j\}_{I\times J}.\cr}
$$
}
\rema{NN}
{
Let $\c C$ be an abelian category,
then $\Pro(\c C)$ is abelian (see [AM, appendix]).
Moreover if $\c C$ has enough injectives
the same is true for $\Pro(\c C)$ (see [J], or [AM]). 
This result mainly concerns the description
of ${\rm Ker}(f)$ and ${\rm Coker}(f)$, where $f$
is a morphism in  $\Pro(\c C)$, done in [AM].
So $\nu(\c O_X)$ is an abelian category. 
Let $\Pro(\QCoh(\c O_X))$ be the Pro-category of
quasi-coherent $\c O_X$-modules.
It is an abelian category and it has enough injectives
(since $\QCoh(\c O_X)$ has enough injectives),
moreover $\nu(\c O_X)$ is
a full thick subcategory of 
$\Pro(\QCoh(\c O_X))$.
In fact $\Coh(\c O_X)$ is a full thick subcategory of 
$\QCoh(\c O_X)$, and it is easy to prove that the same is true for
their Pro-categories ([F2]).
}
\smallskip

\rema{NN}
{
Let denote by $D^+_{\nu(\c O_X)}(\Pro(\QCoh(\c O_X)))$
the derived category of $\Pro(\QCoh(\c O_X))$
with cohomology bounded below and in $\nu(\c O_X)$.
Then 
$D^+_{\nu(\c O_X)}(\Pro(\QCoh(\c O_X)))$
is equivalent to the derived category
$D^+(\nu(\c O_X))$.
}
\smallskip

\smallskip
\section{Stratified Pro-Modules}
In this paper we consider $X$ a smooth algebraic variety over $\Bbb C$.
We will denote by $\{\c P^m_X\}_{\b Z}$ the projective system of
sheaves of principal parts ([EGA IV]),
by $q_m:\c P^m_X\lra \c O_X$ the map induced  by the diagonal embedding
$X\lra X\times X$ and by $q_{m,n}:\c P^m_X\lra \c P^n_X$ ($m\geq n$)
the maps of the projective system $\{\c P^m_X\}_{\b Z}$.
By definition $\c D_{X,m}=\cHom_{\c O_X}(\c P^m_X,\c O_X)$ and 
$\c D_X=\limind_{m\in \Bbb N}\c D_{X,m}$
is the sheaf of differential operators.
We denote by
$\cOmega^\point_X$ the De Rham complex of algebraic differential forms
and by
$\cTheta^{-i}_X\bydef \cHom_{\c O_X}(\cOmega^i_X,\c O_X)$
its dual.
Moreover let $d\bydef d_X$ be the dimension of $X$;
we denote by 
$\omega_X\bydef \cOmega^d_X$ the sheaf of
differential forms of maximum degree.
\smallskip 

\defi{NN}
{ ([BeO; 2.10]).
Let $X$ be a smooth separated noetherian scheme of finite type over $\Bbb C$ and
let $\c F$ be an $\c O_X$-module. 
A  stratification on $\c F$ is a collection 
(one for any $n\in \Bbb N$) of  
$\c P_{X}^n$-linear isomorphisms  
$$
\varepsilon_{\c F,n}:
\c P_{X}^n\otimes_{\c O_X}\c F\lra \c F\otimes_{\c O_X}\c P_{X}^n
$$ 
such that
$\varepsilon_{\c F,n}$ and $\varepsilon_{\c F,m}$ are compatible via $q_{n,m}$ for each
$m\leq n$,  the map
$\varepsilon_{\c F,0}$ is the identity, and the cocycle condition holds. 
}
\smallskip

\prop{NN}
{\rif{eqstrat}([BeO; 2.11]).
Let $\c F$ be an 
$\c O_X$-module, 
the following are equivalent:
\indent
\item{i)}
there is a collection of maps 
$$
s_{\c F,n}:\c F\lra \c F\otimes_{\c O_X}\c P^n_{X}
$$ 
``right'' $\c O_X$-linear such that  
$ s_{\c F,0}=id_{\c F}$,  
$(id_{\c F}\otimes_{\c O_X}q_{m,n})\circ s_{\c F,m}=s_{\c F,n}$ 
and  
$(s_{\c F,n}\otimes_{\c O_X}id_{\c P^m_{X}})\circ s_{\c F,m}= 
(id_{\c F}\otimes_{\c O_X}\delta^{n,m})\circ s_{\c F,m+n}$ 
(see [EGA IV,16.8.9.1] for the definition of $\delta^{m,n}$);
\item{i')}
there is a collection of maps 
$$
s'_{\c F,n}:\c F\lra \c P^n_{X}\otimes_{\c O_X}\c F
$$ 
``left'' $\c O_X$-linear such that  
$ s'_{\c F,0}=id_{\c F}$,  
$(q_{m,n}\otimes_{\c O_X}id_{\c F})\circ s'_{\c F,m}=s'_{\c F,n}$ 
and  
$(id_{\c P^m_{X}}\otimes_{\c O_X}s'_{\c F,n})\circ s'_{\c F,m}= 
(\delta^{n,m}\otimes_{\c O_X}id_{\c F})\circ s'_{\c F,m+n}$ ;
\item{ii)} 
$\c F$ is a stratified module;
\item{iii)}
$\c F$ is a left $\c D_{X}$-module, where 
$\c D_{X}$ is the sheaf of rings of differential operators. 
} 
\proof{ 
$i)\Leftrightarrow ii)$ [BeO, 2.11]. 
\endgraf
$i)+ii)\Rightarrow i')$ and  $i')+ii)\Rightarrow i)$.
\endgraf 
$iii)\Leftrightarrow i)$ 
Let  
$\c D_{X}\otimes_{\c O_X}\c F\buildrel {m_{\c F}}\over\lra \c F$ 
be the multiplication on the $\c D_{X}$-module $\c F$ 
then  
$$ 
\eqalign{ 
m_{\c F}\in {\rm Hom}_{\c O_X} 
(\c D_{X}\otimes_{\c O_X}\c F,\c F)			& 
=\limpro_{m\in \Bbb Z}{\rm Hom}_{\c O_X} 
(\c D_{X,m}\otimes_{\c O_X}\c F, \c F)		\cr 
							& 
=\limpro_{m\in \Bbb Z}{\rm Hom}_{\c O_X} 
(\c F,\c F\otimes_{\c O_X}\c P_{X}^m).\cr} 
$$ 
The associative diagram induces the diagram for co-associativity,  and the identity 
diagram induces that of the co-identity. 
} 
\smallskip

Stratified $\c O_X$-modules form a category which we denote by
$\c O_X\hbox{-}\Strat$.
Morphisms are $\c O_X$-linear maps which respect the stratifications.
We are now interested only in coherent objects so
$\Coh(\c O_X)$-Strat will denote the full subcategory of
$\c O_X$-Strat whose objects are coherent.
\par
We want to extend this category to Pro-objects, in order to obtain
a category dual to that of quasi-coherent 
(so $\Ind(\Coh(\c O_X))$) right $\c D_X$-modules.
\par
The naive way would be that of taking simply the Pro-category
$\Pro(\c O_X\hbox{-}\Strat)$,
but in this way we obtain Pro-objects which have a stratification
at any ``level" while we need a larger category, that of stratified
Pro-objects defined as follow.
\par

 
\defi{NN} 
{ 
Let $ \nu(\c P_{X}^\cdot) $  be the category whose objects 
are Pro-coherent $\c O_X$-modules $\{\c F_h\}_{H}$ endowed with a 
stratification that is a morphism of Pro-objects
$$\xymatrix{ 
\{\c F_h\}_H 
\ar[r]^(0.3){s_{\{\c F_h\}_H}}&\{\c F_h\}_H\otimes_{{\c 
O}_{X}}\{{\c P}_{X}^{m}\}_{\Bbb Z}\cr 
}$$ 
which make the co-identity diagram  
$$\xymatrix{ 
\{\c F_h\}_H 
\ar[r]^(0.3){s_{\{\c F_h\}_H}}	 
\ar[rd]_{id_{\{\c F_h\}_H}}					& 
\{\c F_h\}_H\otimes_{{\c O}_{X}}\{{\c P}_{X}^{m}\}_{\Bbb Z} 
\ar[d]^{id_{\{\c F_h\}_H} 
\otimes \{q_{m}\}_{\Bbb 
 Z} 	}							\cr 
  & \{\c F_h\}_H\cr}  
\leqno{\bf (\NNN)}\rif{coid}
$$ 
and the co-associative one
$$ 
\xymatrix{ 
\{\c F_h\}_H 
\ar[r]^{s_{\{\c F_h\}_H}} 
\ar[d]^{s_{\{\c F_h\}_H}} 					&  
\{\c F_h\}_H\otimes \{{\c P}_{X}^{m}\}_{\Bbb Z} 
\ar[d]^{s_{\{\c F_h\}_H}\otimes id_{\{{\c P}_{X}^{m}\}_{\Bbb Z}}} \cr 
\{\c F_h\}_H \otimes \{{\c P}_{X}^{m}\}_{\Bbb Z} 
\ar[r]_{id_{\{\c F_h\}_H}\otimes s_{\{{\c P}_{X}^{m}\}_{\Bbb Z}}}	&  
\{\c F_h\}_H 
\otimes \{{\c P}_{X}^{m}\}_{\Bbb Z}\otimes\{{\c P}_{X}^{m}\}_{\Bbb Z}\cr 
}\leqno{\bf (\NNN)}\rif{coass}$$ 
commutative. 
(This is simply the category of $\{\c P^m_{X}\}_{\Bbb Z}$-co-modules 
in the category $\nu(\c O_X)$).
\endgraf 
By definition  
$ 
s_{\{{\c P}_{X}^{m}\}_{\Bbb Z}}:=\{\delta^{m,n}\}_{\b Z\times \b Z}
$  
is the map inducing the stratification on  the ``right"
$\{{\c P}_{X}^{m}\}_{\Bbb Z}=p_1^\ast(\c O_X)$ 
(see [BeO, Remark 2.13], [G]). 
A morphism of Pro-objects 
$f:\{\c F_h\}_H\lra \{\c G_k\}_K$  
is a morphism in $\nu(\c P_{X}^\cdot)$  
if and only if the diagram 
$$\xymatrix{ 
\{\c F_h\}_H 
\ar[r]^{f}	 
\ar[d]^{s_{\{\c F_h\}_H}}				& 
\{\c G_k\}_K						 
\ar[d]^{s_{\{\c G_k\}_K}}				\cr 
\{{\c P}_{X}^{m}\}_{\Bbb Z}\otimes \{\c F_h\}_H 
\ar[r]_{id
\otimes f}	& 
\{{\c P}_{X}^{m}\}_{\Bbb Z}\otimes \{\c G_k\}_K		\cr 
}\leqno{\bf (\NNN)}\rif{comor}$$ 
commutes.
We denote by 
$\cHom_{\c P_{X}^\cdot}(\{\c F_h\}_H,\{\c G_k\}_K)$
(or $\cHom_{\Strat}(\{\c F_h\}_H,\{\c G_k\}_K)$
the sheafified version of the set of morphisms in 
$\nu({\c P_{X}^\cdot})$.
\endgraf
We denote by
$C^\ast(\c P^\cdot_X)$ (resp.
$K^\ast(\c P^\cdot_X)$, resp.
$D^\ast(\c P^\cdot_X)$) with $\ast\in \{+,-,b\}$ 
the category of  complexes (bounded below, bounded above, bounded)
(resp. up to homotopy, resp. up to quasi-isomorphisms) in $\nu(\c O_X)$.
} 
\smallskip

\rema{NN}
{\rif{localbase}[BeO; 2.2, 2.3].
In our setting $X$ is smooth, hence the sheaves $\{\c P^m_X\}_{\b Z}$ are locally free 
of finite type.
A base of $\c P^m_X$ (for both left and right $\c O_X$-module structures)
is 
$\{\xi_1^{\alpha_1}\cdots\xi_d^{\alpha_d}| \alpha_1+\cdots +\alpha_d\leq n; \;\alpha_i\in \Bbb N\}$
where
$\xi_i\bydef 1\otimes x_i-x_i\otimes 1$ and $\b I= 1\otimes 1$ in local coordinates.
By this description the ``right" stratification on
the Pro-system  $\{\c P^m_X\}_{\b Z}$ is given  by the morphisms
$$
\matrix{
\delta^{m,p}: \c P^m_X &   \lra &  \c P^{m-p}_X\otimes_{\c O_X} \c P^p_X \hfill&\cr
	\hfill	\b I		& \longmapsto & \b I\otimes \b I\hfill& \cr
\hfill\xi_i &	\longmapsto & \b I\otimes \xi_i +\xi_i\otimes \b I\hfill& \cr
\hfill \xi_1^{\alpha_1}\cdots\xi_d^{\alpha_d}&\longmapsto &
\delta^{m,p}(\xi_1)^{\alpha_1}\cdots \delta^{m,p}(\xi_d)^{\alpha_d}&
\forall \alpha_1,\dots,\alpha_d\in\b N,\; \alpha_1+\cdots +\alpha_d\leq m\cr
}
$$
Moreover the sheaf $\cOmega^1_X$ is simply the sub-sheaf 
of $\c P^1_X$ generated by $\xi_i$ for $i=1,\dots,d$.
}
\smallskip
 
\defi{NN}
{
A stratified Pro-module is induced if it is isomorphic to
$\{\c P^m_X\}_{\b Z}\otimes_{\c O_X}\c L$,
for some $\c L\in\nu(\c O_X)$, endowed
with the stratification induced by  the canonical one on
$\{\c P^m_X\}_{\b Z}$ (see [G, 6.3]).
We denote by $\nu_i(\c P^\cdot_X)$ the 
full subcategory of $\nu(\c P^\cdot_X)$
whose objects are induced.
We denote by
$C^b_i(\c P^\cdot_X)$ (resp.
$K^b_i(\c P^\cdot_X)$, resp.
$D^b_i(\c P^\cdot_X)$) the category of bounded complexes 
(resp. up to homotopy, resp. up to quasi-isomorphisms) in $\nu_i(\c O_X)$.
}
\smallskip
 
 
\prop{NN} 
{ \rif{coext}
The category
$\nu(\c P^\cdot_X)$ is an abelian category,
small filtering 
projective limits are representable  and exact. 
The forgetful functor 
$$for: \nu(\c P^\cdot_X)\lra \nu(\c O_X)$$ 
has a right adjoint  
$$\matrix{ 
\Q^0_X\bydef \{\c P^m_{X}\}_{\b Z}\otimes_{\c O_X}\_:&
\nu(\c O_X) &\lra & \nu(\c P^\cdot_X)\hfill\cr 
&\hfill\{\c F_h\}_H&\longmapsto&\{\c P_{X}^m\}_{\Bbb Z}\otimes_{\c O_X}\{\c 
F_h\}_H\cr}$$
which takes image into $\nu_i(\c P^\cdot_X)$.
}  
\proof{ 
${\rm Kernels}$ and ${\rm cokernels}$ in $\nu(\c P_{X}^\cdot)$ 
are those of $\nu(\c O_X)$
endowed with the induced stratification and
for any morphism $f$ in $\nu(\c P_{X}^\cdot)$,
the image of $f$ is isomorphic to its co-image.
So  $\nu(\c P_{X}^\cdot)$ 
is an abelian category and the forgetful functor is exact.
Small filtering limits are representable and exact
because they are representable in
$\nu(\c O_X)$ and they  have canonical stratifications.
\endgraf
The map 
$$  
\matrix{ 
{\c Hom}_{ \Strat} 
(\{ \c F_h\}_H,\{ \c P_{X}^m\}_{\b Z}\otimes_{\c O_X}\{\c G_k\}_K) & 
\buildrel{\alpha}\over{\ra} & 
\cHom_{\nu(\c O_X)}(for(\{\c F_h\}_H),\{\c G_k\}_K)\cr 
\hfill f & \mapsto & (\{q_m\}_{\b Z}\otimes_{\c O_X}id_{\{\c G_k\}_K})\circ f\hfill\cr 
} 
$$ 
(co-extension of scalars)
is a bijection whose inverse is the map
$$ 
\matrix{ 
\cHom_{\nu(\c O_X)}(for(\{\c F_h\}_H),\{\c G_k\}_K)			& 
\buildrel{\beta}\over{\ra} 								& 
{\c Hom}_{ \{ \c P_{X}^m \}_{\b Z}} 
(\{ \c F_h\}_H,
\{ \c P_{X}^m\}_{\b Z}\otimes_{\c O_X}\{\c G_k\}_K) 		\cr 
\hfill g											&
\mapsto 											& 
(id_{\{\c P_{X}^m\}_{\b Z}}\otimes g)\circ 
\{s'_{\c F,m}\}_{\b Z}.\hfill	\cr 
} 
$$ 
Clearly $\alpha(f)$ is a morphism in $\nu(\c O_X)$;
on the other hand in order to prove that $\beta(g)$ respects
the stratifications 
it is sufficient to 
remark that $\Q^0_X$ is a functor so
the map $(id_{\{\c P_{X}^m\}_{\b Z}}\otimes g)$ respects the
canonical stratifications on $\{\c P_{X}^m\}_{\b Z}\otimes _{\c O_X}\_$.
} 
\smallskip

\rema{NN}
{
Proposition \cite{coext} also  holds true on replacing 
$\nu(\c O_X)$ by the category of Pro-quasi-coherent
$\c O_X$-modules $\Pro(\mu(\c O_X))$,
and
$\nu(\c P^\cdot_X)$ by the category of
stratified Pro-quasi-coherent
$\c O_X$-modules denoted $\Pro(\mu(\c P^\cdot _X))$.
}
\smallskip
 
\coro{NN}
{
Any object in  $\Pro(\mu(\c P^\cdot _X))$
induced by an injective object in $\Pro(\mu(\c O_X))$
is injective. Moreover $\Pro(\mu(\c P^\cdot _X))$
has enough injectives.
}
\proof{
Let $\c E$ be an injective Pro-quasi-coherent 
$\c O_X$-module, then the functor
$$
\cHom_{\Strat}(\_,\{\c P^m_X\}\otimes_{\c O_X}\c E)\cong
\cHom_{\Pro(\mu(\c O_X))}(for(\_),\c E)
$$
is exact because $for(\_): \Pro(\mu(\c P^\cdot _X))\ra \Pro(\mu(\c O_X))$
is  exact and
$\c E$ is injective.
\endgraf
For each $\c N\in \Pro(\mu(\c P^\cdot _X))$, there exists
$\c I\in\Pro(\mu(\c O_X))$ and an injective map
$i: for(\c N)\hookrightarrow \c I$ in $ \Pro(\mu(\c O_X))$.
Then the map 
$\beta(i): \c N \lra \{ \c P_{X}^m\}_{\b Z}\otimes_{\c O_X}\c I$
is an injective map in  $\Pro(\mu(\c P^\cdot _X))$
(and $\{ \c P_{X}^m\}_{\b Z}\otimes_{\c O_X}\c I$
is injective in $\Pro(\mu(\c P^\cdot _X))$.}
\smallskip


\coro{NN}
{
{\bf{Derived co-extension of scalars.}}\endgraf
 Let $\c F\in\nu(\c P^\cdot _X)$ and
$\c G\in\nu(\c O_X)$:
$$ 
\bR\cHom_{\Strat} 
(\c F,\{\c P_{X}^m\}_{\b Z}\otimes_{\c O_X}\c G)  
\buildrel{\cong}\over \lra \bR{\c Hom}_{\nu(\c O_X)}(for(\c F),\c G); 
$$ 
is a quasi-isomorphism.
}
\proof{
Let denote by $E^\point(\c G)$ an injective resolution of $\c G$
in $\Pro(\mu(\c O_X))$.
Then
$$
\eqalign{
\bR\cHom_{\Strat}(\c F,\{\c P_{X}^m\}_{\b Z}\otimes_{\c O_X}\c G)&=
\cHom_{\Strat}(\c F,\{\c P_{X}^m\}_{\b Z}\otimes_{\c O_X}E^\point(\c G))\cong\cr
&\cong \cHom_{\Pro(\mu(\c O_X))}(for(\c F),E^\point(\c G))=\cr
& = \bR\cHom_{\nu(\c O_X)}(for(\c F),\c G).\cr}
$$
}
\smallskip
  
Let consider the De Rham functor
$$\bR\cHom^\cdot_{\Strat}(\c O_X,\_):D^b(\nu(\c P^\cdot_X))\lra D^b(\Bbb C_X)
$$
where $\Bbb C_X$ denotes the category of sheaves in $\Bbb C$-vector spaces.
Then for any
$\c M^\cdot\in D^b(\nu(\c P^\cdot_X))$ deriving it in $\Pro(\mu(\c P^\cdot_X))$
we have:
$$
\eqalign{
\bR\cHom^\cdot_{\Strat}(\c O_X,\c M^\cdot) \cong &
\bR\cHom^\cdot_{\Strat}(\c D_X\otimes \cTheta^\cdot_X,\c M^\cdot) \cong\cr
\cong & \cHom^\cdot_{\Strat}(\c D_X\otimes \cTheta^\cdot_X,\c M^\cdot) \cong\cr
\cong & \cOmega^\cdot_X\otimes^\cdot \c M^\cdot\cr
}
$$
The complex $ \cOmega^\cdot_X\otimes^\cdot \c M^\cdot$
is a complex of Pro-coherent-$\c O_X$-modules but its differentials are not
$\c O_X$-linear.
\endgraf
In the following we will define the category $\nu(\c O_X)\hbox{-}\Diff_X$
wherein the functor
$\cHom^\cdot_{\Strat}(\c O_X,\_)$ has its image in a fully faithful way.
So the De Rham functor will have its image
in a suitable localization of  $\nu(\c O_X)\hbox{-}\Diff_X$.
\smallskip

\theo{NN}
{
Induced stratified Pro-modules are acyclic for the functor 
$\cHom_{\Strat}(\c O_X,\_)$.
For $\c M$, $\c N$ such modules 
$$
\c M^{\nabla}\bydef 
\bR\cHom_{\Strat}(\c O_X,\c M)=
\cHom_{\Strat}(\c O_X,\c M)$$
and the morphism
$$
\cHom_{\Strat}(\c M,\c N)\lra \cHom_{\Bbb C_X}(\c M^{\nabla},\c N^{\nabla})
\leqno{\bf (\NNN)}\rif{injmor}
$$
is injective.
}
\proof{
By hypothesis there exist $\c L$, $\c L'$ in $\nu(\c O_X)$ such that 
$\c M\cong\{\c P^m_X\}_{\b Z}\otimes_{\c O_X}\c L$
and
$\c N\cong\{\c P^m_X\}_{\b Z}\otimes_{\c O_X}\c L'$.
Then
$$\eqalign{
\bR\cHom_{\Strat}(\c O_X,\c M)=&
\bR\cHom_{\Strat}(\c O_X,\{\c P^m_X\}_{\b Z}\otimes_{\c O_X}\c L)\cong\cr
\cong& \bR\cHom_{\nu({\c O_X})}(\c O_X,\c L)\cong\cr
\cong& \cHom_{\nu({\c O_X})}(\c O_X,\c L)\cong\cr
\cong& \cHom_{\Strat}(\c O_X,\{\c P^m_X\}_{\b Z}\otimes_{\c O_X}\c L)=\cr
=&\cHom_{\Strat}(\c O_X,\c M)\cr}
$$
which proves the first assertion.
\endgraf
For the second statement let consider the map
$$\eqalign{
\cHom_{\Strat}(\c M,\c N)=&
\cHom_{\Strat}(\{\c P^m_X\}_{\b Z}\otimes_{\c O_X}\c L,
\{\c P^m_X\}_{\b Z}\otimes_{\c O_X}\c L')\cong\cr
\cong&\cHom_{\nu({\c O_X})}(\{\c P^m_X\}_{\b Z}\otimes_{\c O_X}\c L,\c L')
\lra \cHom_{\Bbb C_X}(\c L,\c L')\cr
}\leqno{\bf (\NNN)}\rif{injmorpro}
$$
obtained by composition with the stratification morphism
$s'_{\c L}:\c L\ra \{\c P^m_X\}_{\b Z}\otimes_{\c O_X}\c L$.
It is injective because the image of $s'_{\c L}$ generates 
$\{\c P^m_X\}_{\b Z}\otimes_{\c O_X}\c L$ as 
Pro-coherent-$\c O_X$-module.\endgraf
We note that this theorem is the analogue (for Pro-objects)
of Saito's Lemma [S.2; 1.2].
}
\smallskip

\section{Differential Complexes of Pro-modules}
\defi{NN}
{
Let $\{\c L_i\}_I$ and $\{\c L'_j\}_J$ be two Pro-coherent $\c O_X$-modules.
  The sheaf of differential operators, which we denoted by
$\cHom_{\Diff}(\{\c L_i\}_I,\{\c L'_j\}_J)$, is
the image of the injective map
(\cite{injmorpro}).
So
$$\eqalign{
\cHom_{\Diff_X}(\{\c L_i\}_I,\{\c L'_j\}_J)\bydef&
\cHom_{\Strat}(\{\c P^m_X\}_{\b Z}\otimes_{\c O_X}\c L,
\{\c P^m_X\}_{\b Z}\otimes_{\c O_X}\c L')\cong\cr
\cong&
\cHom_{\nu(\c O_X)}(\{\c P^m_X\}_{\b Z}\otimes_{\c O_X}\{\c L_i\}_I,\{\c L'_j\}_J)\bydef\cr
\bydef&\limpro_J\limind_I\limind_\Bbb Z\cHom_{\c O_X}(\c P^m_X\otimes_{\c O_X}\c L_i,
\c L'_j)\cong\cr
\cong&\limpro_J\limind_I\cHom_{\Diff_X}(\c L_i,\c L'_j).\cr
}
$$
We recall that for $\c F\in\Coh(\c O_X)$ and $\c G$ an $\c O_X$-module,
the sheaf 
$\cHom_{\Diff_X}(\c F,\c G)$ is isomorphic to 
${\limin_{\Bbb Z}}\cHom_{\c O_X}(\c P^m_X\otimes_{\c O_X}\c F,\c G)\cong
{\limin_{\Bbb Z}}\cHom_{\c O_X}(\c F,\c G\otimes_{\c O_X}\c D_{X,m})$.
\endgraf
We denote by $\nu(\c O_X)\hbox{-}\Diff_X$ the additive category 
whose objects are Pro-coherent-$\c O_X$-modules and whose morphisms are 
differential operators (sometimes called differential complexes).
We have a  functor
$$
\matrix{
{\Q}^0_X:&\nu(\c O_X)\hbox{-}\Diff_X & \lra &\nu_i(\c P^\cdot_X)\hfill\cr
		&\hfill\c L		&\longmapsto&\{\c P^m_X\}_{\b Z}\otimes_{\c O_X}\c L\cr
}
$$
which extends that of Proposition \cite{coext}.
This functor was firstly introduced by A. Grothendieck in [G; 6.2] and it is called the
formalization functor or linearization. 
By (\cite{injmor}) it is an equivalence of categories.
}\smallskip

\rema{NN}{
If we restrict the formalization functor to differential
complexes $\c L$ whose objects are coherent $\c O_X$-modules,
then the Pro-objects $Q^0_X(\c L)$ are always of Artin-Rees type.
Moreover any morphism of Pro-obejcts between two such objects
is  necessarily of Artin-Rees type (see [G;6.2]).
}\smallskip

\defi{NN}{
Let $C^b(\nu(\c O_X),\Diff_X)$ be the category of bounded
complexes in $\nu(\c O_X)\hbox{-}\Diff_X$.
Let
$D^b(\nu(\c O_X),\Diff_X)$ be the category obtained from 
$C^b(\nu(\c O_X),\Diff_X)$ by inverting 
${\Q}^0_X$-quasi-iso\-morphisms.
This is a triangulated category with the usual shift functor 
and distinguished triangles those induced by the usual mapping cones.
We remark that this localizing procedure was first introduced 
in [AB, Appendix C] following an idea of P. Berthelot.
\endgraf
We obtain a localized equivalence of categories
$$
\Q^0_X: D^b(\nu(\c O_X),\Diff_X)\lra D^b_i(\c P^\cdot_X).
$$
This functor would be the ``dual" of Saito $\tDR^{-1}_{X}$
functor.
} 
\smallskip
 
 \rema{NN}
 {
 The morphism (\cite{injmor}) is induced by the following commutative
 diagrams which are adjoint to those of $\c D_X$-modules in the smooth
 case ([S.2; (1.4.1)]):
 $$
 \xymatrix{
 \c L \ar[r]^{P }\ar[d]_{d^1\otimes id_{\c L}} & \c L'	\ar[d]^{d^1\otimes id_{\c L'}} 
 &\;	& 
 \c L\ar[r]^{P} \ar[d] _{d^1\otimes id_{\c L}} & \c L'\cr
\{\c P^m_X\}_{\b Z}\otimes_{\c O_X}\c L
\ar[r]_{\Q^0_X(P)} 
&\{\c P^m_X\}_{\b Z}\otimes_{\c O_X}\c L'
&	\;&
\{\c P^m_X\}_{\b Z}\otimes_{\c O_X}\c L\ar[ur]_{\Q^0_X(P)} & \cr
 }\leqno{\bf (\NNN)}\rif{comdiadiop}
 $$
 for any $\Q^0_X(P)\in \Hom_{\Strat}({\Q}^0_X(\c L),{\Q}^0_X(\c L'))
\cong \Hom_{\nu({\c O_X})}({\Q}^0_X(\c L),\c L')$.
The map $d_1:\c O_X\lra \c P^m_{X}$ is that induced by
the second projection of $X\times X$ into $X$.
 }
 \smallskip
 
\defi{NN} 
{ 
An object in $D^b(\nu(\c O_X),\Diff_{X})$ is said to be perfect if it is
locally isomorphic to a bounded complex whose elements are 
locally free $\c O_X$-modules of finite rank.
We denote by $D^b_p(\nu(\c O_X),\Diff_{X})$
the category of bounded perfect complexes in $D^b(\nu(\c O_X),\Diff_{X})$.
Then any object in $D^b_p(\nu(\c O_X),\Diff_{X})$
may be represented as an object in $C^b(\Coh(\c O_X),\Diff_{X})$
(see [L] for definition of perfect objects).
}
\smallskip

\defi{NN}{\rif{HLCX}
{\bf Herrera-Lieberman differential complexes.}\endgraf
([HL, \S 2] or [Be, II.5]).
Let  $C^b_{1}(\nu(\c O_X),\Diff_{X})$ denote the category 
of bounded complexes of differential operators of order 
at most one,
that is:
\endgraf
	\item{i)} the objects of $C^b_{1}(\nu(\c O_X),\Diff_{X})$ 
	are complexes whose terms are Pro-coherent-$\c O_{X}$-modules and 
	whose differentials are 			
	differential operators of order at most one;
\endgraf
	\item{ii)} morphisms between such complexes 
	are morphisms of complexes
	which are $\c O_{X}$-linear maps.
\endgraf
The category 
$C^b_{1}(\Coh(\c O_X),\Diff_{X})$ 
is the full subcategory 
of $C^b_{1}(\nu(\c O_X),\Diff_{X})$
whose objects are coherent modules.	
\endgraf
We denote by $D^b_{1}(\nu(\c O_X),\Diff_{X})$ 	
the category obtained form $C^b_{1}(\nu(\c O_X),\Diff_{X})$
by inverting ${\Q}^0_X$-quasi-isomorphisms
in
$C^b_{1}(\nu(\c O_X),\Diff_{X})$.
Thus we have the functors
$$
\eqalign{
\lambda: D^b_{1}(\nu(\c O_X),\Diff_{X})&\lra D^b(\nu(\c O_X),\Diff_{X})\cr
{\Q}^0_{X,1}: D^b_{1}(\nu(\c O_X),\Diff_{X}) &\lra D^b(\c P^\cdot_X)\cr
}
$$
where ${\Q}^0_{X,1}\bydef {\Q}^0_X\opp \lambda$.
}
\smallskip

\section {De Rham Functor} 

\defi{NN}
{
Let $\c M\in \nu(\c P^\cdot_X)$ and
$$
\ov{\DR}_X(\c M)\bydef 
\cOmega^\point_X\otimes_{\c O_X} \c M=
$$
$$=
\xymatrix{
[0\ar[r] &
\mathop{\c M}\limits^{0}_{\phantom{0}} \ar[r]  &
\mathop{\cOmega^1_X\otimes_{\c O_X}\c M}\limits^{1}_{\phantom{1}}  \ar[r] &
\cdots \ar[r]  
&
\mathop{\cOmega^d_X\otimes_{\c O_X}\c M}\limits^{d}_{\phantom{d}} 
\ar[r]&
0].  \cr}
\leqno{\bf (\NNN)}\rif{ovDR}
$$
The differentials are defined using the stratification map $s'_{\c M}$ (see Proposition
\cite{eqstrat}) and the projection $\c P^1_X\lra \cOmega^1_X$.
The complex $\ov{\DR}_X(\c M)$
belongs to $C^b(\nu(\c O_X),\Diff_X)$
and in particular it is also an object of  $C^b_1(\nu(\c O_X),\Diff_X)$. 
We define the functors
$$
\matrix{
\ov{\DR}_{1,X}: &C^b(\c P^\cdot_X)&\lra&
C^b_1(\nu(\c O_X),\Diff_X)\hfill\cr
& \hfill \c M^\point &\longmapsto &
(\cOmega^\point_X\otimes_{\c O_X}\c M^{\point})_{\ss {tot}}=:
\cOmega^\point_X\otimes^{\point}_{\c O_X}\c M^{\point}.}
$$
and $\ov{\DR}_X=\lambda_1\opp \ov{\DR}_{1,X}$.
}
\smallskip

We want to prove that this De Rham functor sends
the multiplicative system of qis in $C^b(\c P^\cdot_X)$ into
the multiplicative system of $\Q^0_X$-qis in $C^b(\nu(\c O_X),\Diff_X)$.
In order to prove this result we need the following version of
the crystalline Poincar\'e  lemma
$\c O_X\buildrel{\rm qis}\over\lra {\Q}^0_X \ov{\DR}_X(\c O_X)$.
This lemma may be found in  [G, 6.5] and in 
[BeO, 6.12] where Berthelot Ogus proved a filtered version.
We give here a simple proof of the result we need.
We note that our proof also works well in characteristic p
using the formalism of divided powers.
\smallskip
Let us remark:

\rema{NN}
{[Be.O, 2.13].
The Pro-object
$\{\c P^m_X\}_{\b Z}$ admits two different stratifications depending on
the $\c O_X$-module structure we chose on it.
We consider on $\{\c P^m_X\}_{\b Z}$
its ``left" $\c O_X$-structure, (that given by $p_0$), (so
its ``right" $\c O_X$-structure may be used in the tensor product with the
De Rham complex).
This is the construction of  Grothendieck linearization.
In this case $\{\c P^m_X\}_{\b Z}$ is endowed with the stratification 
$\theta:\c P^{m+n}_X\ra {\c P^n_X}\otimes_{\c O_X}{\c P^m_X}$
sending $(f\otimes g)\mapsto 1\otimes g\otimes 1\otimes f$. 
\endgraf
On the other hand, if we consider the ``right'' structure on $\{\c P^m_X\}_{\b Z}$, 
the stratification is given by the map
$\delta$.
}\smallskip

\lemm{NN}{\rif{crispo}
The linearized De Rham complex is a resolution of $\c O_X$
$$
\c O_X \buildrel{d^0}\over\lra {\Q}^0_X \ov{\DR}_X(\c O_X) 
$$ 
in 
$C^b(\c P^\cdot_X)$.
In fact the complex
$$\xymatrix{
\c O_X\ar[r]^(0.4){d_0}& \{\c P^m_X\}_{\b Z} \ar[r]^(0.4){\ov{\nabla}^0}&\{\c P^m_X\}_{\b Z}\otimes_{\c O_X}\cOmega^1_X\ar[r]^(0.7){\ov{\nabla}^1}&
\cdots \ar[r]^(0.3){\ov{\nabla}^{d-1}} & \{\c P^m_X\}_{\b Z}\otimes_{\c O_X}\cOmega^d_X\cr
}\leqno{\bf (\NNN)}\rif{Poincris}
$$
is exact and
thus locally homotopic to zero since its terms are locally free.
}
\proof{
First of all we remark that ${\Q}^0_X \ov{\DR}_X(\c O_X)$
is a complex in $C^b(\c P^\cdot_X)$ and
the map
$\c O_X \buildrel{d^0}\over\lra \{\c P^m_X\}_{\b Z}$
respects the stratifications (see the remark given below regarding the stratification on $\{\c P^m_X\}_{\b Z}$).
It is evident that $\ov{\nabla}^0\circ d^0=0$.

The complex (\cite{Poincris}) as a complex of $\nu(\c P^\cdot_X)$  is represented
by
$$
\xymatrix{
\{\c O_X
\ar[r]^{d_0}& 
\c P^n_X \ar[r]^(0.4){\ov{\nabla}^0_{(n)}}&
\c P^{n-1}_X\otimes_{\c O_X}\cOmega^1_X\ar[r]^(0.7){\ov{\nabla}^1_{(n)}}&
\cdots \ar[r]^(0.3){\ov{\nabla}^{d-1}_{(n)}} &
\c P^{n-d}_X\otimes_{\c O_X}\cOmega^d_X\}_{n\in \Bbb N}.\cr
}\leqno{\bf {(\NNN .n)_{n\in \Bbb Z}}}\rif{ProPoincris}
$$
We will prove by induction that (\cite{ProPoincris}.n) is exact for each $n\in \Bbb N$.
\endgraf
First of all we render the $\c O_X$-linear differentials on the complex explicit
by the use of the basis given by $\xi_i$ (see remark \cite{localbase}).
Then $d^0$ is the map
$$
\matrix{
d^0: \c O_X &\lra & \c P^n_X\hfill\cr
\hfill 1&\longmapsto& \Bbb I\bydef 1\otimes 1\hfill\cr
\hfill f&\longmapsto& f\otimes 1=f\b I\hfill\cr
}
$$
while $\ov{\nabla}\bydef \Q^0_X(\nabla)$ is  the linearization 
of the De Rham complex 
($
\xymatrix{
\c O_X\ar[r]^{\nabla^0}&\cOmega^1_X\ar[r]^{\nabla^1}&\cdots\ar[r]^{\nabla^{d-1}}&\cOmega^d_X\cr}
$)
obtained as
$$
\xymatrix{
\c P^n_X\otimes_{\c O_X}\cOmega^p_X\ar[r]^{\ov{\nabla}^p_{(n)}}
\ar[d]_{\delta^{n,1}\otimes{id_{\cOmega^p_X}}}&
\c P^{n-1}_X\otimes_{\c O_X}\cOmega^{p+1}_X\cr
\c P^{n-1}_X\otimes_{\c O_X}\c P^{1}_X\otimes_{\c O_X}\cOmega^p_X
\ar[ru]_{id_{\c P^{n-1}_X}\otimes\ov{\nabla}^p_{(1)}}&\cr
}
$$
where $\ov{\nabla}^p_{(1)}$ is
$$\matrix{
\hfill \ov{\nabla}^p_{(1)}: \c P^1_X\otimes_{\c O_X}\cOmega^p_X & \lra &\cOmega^{p+1}_X\hfill& \cr
\hfill (f\otimes g)\otimes \omega &\longmapsto &f\nabla^p(g\omega) 
\hfill &\rm {so\; that\; in\; local\; coordinates\; we\; have}
\cr
 \hfill \b I\otimes( \xi_{i_1}\wedge\cdots \wedge  \xi_{i_p})&\longmapsto & 0\hfill&\cr
\hfill  \xi_i\otimes(\xi_{i_1}\wedge\cdots \wedge  \xi_{i_p})&\longmapsto &\xi_i\wedge 
\xi_{i_1}\wedge\cdots \wedge  \xi_{i_p}\hfill&\forall i\in\{1,\dots d\}.\;
\cr
}$$
Now we proceed by induction on the ``level" $n$ in order to prove that (\cite{Poincris})
is exact.
For $n=0$ the complex (\cite{ProPoincris}.0) reduces to
$$ 0\lra \c O_X \buildrel{id_{\c O_X}}\over\lra \c O_X \lra 0\lra 0$$
which is obviously exact.
When $n=1$ the complex (\cite{ProPoincris}.1) is
$$ 
0\lra \c O_X \buildrel{d^0}\over\lra \c P^1_X\lra \cOmega^1_X\lra 0
$$
which is  homotopic to zero via the $\c O_X$-linear
homotopism
$$ \matrix{
\c P^1_X&   \buildrel{q_1}\over\lra &\c O_X & &\cOmega^1_X & \lra & \c P^1_X \cr
f\otimes g&\longmapsto&fg& & \xi_i&\longmapsto &\xi_i\cr
}$$
then $\c P^1_X\cong \c O_X\oplus \cOmega^1_X$.
Let us suppose that the complex (\cite{ProPoincris}.n-1) is exact with $\hbox{n}\geq 1$.
Then we consider the diagram
$$
\xymatrix{
& 0\ar[d] &0\ar[d]  &0\ar[d]  &0\ar[d] &0\ar[d]  & \cr
0\ar[r] & 0\ar[d]\ar[r] &\c I^n/\c I^{n+1}\ar[d]\ar[r] &
\c I^{n-1}/\c I^{n}\otimes_{\c O_X}\cOmega^1_X\ar[d]\ar[r] &
\cdots \ar[r]\ar[d] &\c I^{n-d}/\c I^{n-d+1}\otimes_{\c O_X}\cOmega^d_X\ar[d]\ar[r] &0\cr
0\ar[r] & \c O_X\ar[r]\ar[d] &\c P^{n}_X\ar[r]\ar[d] &
\c P^{n-1}_X\otimes_{\c O_X}\cOmega^1_X\ar[r]\ar[d] &\cdots\ar[r]\ar[d] &
\c P^{n-d}_X\otimes_{\c O_X}\cOmega^d_X\ar[r]\ar[d]& 0\cr
0\ar[r] & \c O_X\ar[r]\ar[d] &\c P^{n-1}_X\ar[r]\ar[d] &
\c P^{n-2}_X\otimes_{\c O_X}\cOmega^1_X\ar[r]\ar[d] &\cdots\ar[r]\ar[d] &
\c P^{n-d-1}_X\otimes_{\c O_X}\cOmega^d_X\ar[r]\ar[d]& 0\cr
& 0 &0  &0 &0&0& \cr
}
$$
whose columns are exact. By inductive hypothesis the third row is exact.
Then the second row is exact if and only if the first is.
So we will prove that the complex
$$
0\lra\c I^n/\c I^{n+1}\buildrel{D^0}\over\lra \c I^{n-1}/\c I^{n}\otimes_{\c O_X}\cOmega^1_X
\buildrel{D^1}\over\lra
\cdots\buildrel{D^{d-1}}\over\lra \c I^{n-d}/\c I^{n-d+1}\otimes_{\c O_X}\cOmega^d_X\lra 0
\leqno{\bf {(\NNN .n)}}\rif{GraPoincris}
$$
is exact proving that its identity is homotopic to zero.
\endgraf
We have to construct $\c O_X$-linear maps 
$$
s_p:\c I^{n-p}/\c I^{n-p+1}\otimes_{\c O_X}\cOmega^p_X\lra
\c I^{n-p+1}/\c I^{n-p}\otimes_{\c O_X}\cOmega^{p-1}_X
$$
for any $p=1,\dots, d$ such that in the diagram
$$
\xymatrix{
0\ar[r] &\c I^n/\c I^{n+1}\ar[r]\ar[d]_{id} &
\c I^{n-1}/\c I^n\otimes_{\c O_X}\cOmega^1_X\ar[dl]_{s_1}\ar[r]\ar[d]_{id} &
\cdots\ar[dl]_{s_2}\ar[r]\ar[d]_{id}& 
\c I^{n-d}/\c I^{n-d+1}\otimes_{\c O_X}\cOmega^d_X\ar[dl]_{s_d}\ar[r]\ar[d]_{id} &0\cr
0\ar[r] &\c I^n/\c I^{n+1}\ar[r]&\c I^{n-1}/\c I^n\otimes_{\c O_X}\cOmega^1_X\ar[r]&
\cdots\ar[r]&\c I^{n-d}/\c I^{n-d+1}\otimes_{\c O_X}\cOmega^d_X\ar[r]& 0
}
$$
the identity of $\c I^{n-p}/\c I^{n-p+1}\otimes_{\c O_X}\cOmega^p_X$ would be
$id=D^{p-1}\opp s_p+ s_{p+1}\opp D^p$.
\endgraf\noindent
First we explicitly  write the action of the differentials $D^p$ on a basis: 
$$
\matrix{
\c I^{n-p}/\c I^{n-p+1}\otimes_{\c O_X}\cOmega^p_X & \buildrel{D^p}\over\lra &
\c I^{n-p-1}/\c I^{n-p}\otimes_{\c O_X}\cOmega^{p+1}_X\hfill\cr
\xi_1^{\alpha_1}\cdots\xi_d^{\alpha_d}\otimes \xi_{i_1}\wedge\cdots\wedge\xi_{i_p}&
\longmapsto &
\sum_{j=1}^{d}\alpha_j\xi_1^{\alpha_1}\cdots\xi_j^{\alpha_j-1}\cdots\xi_d^{\alpha_d}\otimes
\xi_j\wedge\xi_{i_1}\wedge\cdots\wedge\xi_{i_p}\cr
}
$$
with $\alpha_1+\cdots+\alpha_d=n-p$.
\endgraf
We note that the map
$$\matrix{
\overbrace{\cOmega^1_X\otimes_{\c O_X}\cdots \otimes_{\c O_X}\cOmega^1_X}^{\hbox{p times}}&\lra &
\overbrace{\cOmega^1_X\otimes_{\c O_X}\cdots \otimes_{\c O_X}\cOmega^1_X}^{\hbox{p times}}\hfill\cr
\alpha_1\otimes\cdots\otimes\alpha_p&\longmapsto & 
\sum_{\sigma\in\Sigma_p}(-1)^{sgn(\sigma)}\alpha_{\sigma(1)}\otimes\cdots\otimes\alpha_{\sigma(p)}\cr
}
$$
induces a map
$\sigma^p:\cOmega^p_X\lra \cOmega^1_X\otimes_{\c O_X}\cOmega^{p-1}_X$.
We define $s_p$ (up to the factor $n(p-1)!$) as the composition
$$
\xymatrix{
\c I^{n-p}/\c I^{n-p+1}\otimes_{\c O_X}\cOmega^p_X
\ar[r]^{n(p-1)! s_p}\ar[d]_{id\otimes \sigma^p}& \c I^{n-p+1}/\c I^{n-p+2}\otimes_{\c O_X}\cOmega^{p-1}_X\cr
\c I^{n-p}/\c I^{n-p+1}\otimes_{\c O_X}\c I/\c I^{2}\otimes_{\c O_X}\cOmega^{p-1}_X
\ar[ru]_{m\otimes id}&\cr
}
$$
where $m$ is the map
$m:\c I^{n-p}/\c I^{n-p+1}\otimes_{\c O_X}\c I/\c I^{2} \lra \c I^{n-p+1}/\c I^{n-p+2}$
induced by the multiplication.
It is well defined because $\c I^{n-p+1}\c I=\c I^{n-p+2}=\c I^{n-p}\c I^{2}$.
We now explicitly calculate $s_p$ on a local basis 
$$
\matrix{
\c I^{n-p}/\c I^{n-p+1}\otimes_{\c O_X}\cOmega^p_X & \buildrel{s_p}\over\lra &
\c I^{n-p+1}/\c I^{n-p+2}\otimes_{\c O_X}\cOmega^{p-1}_X\hfill\cr
\xi_1^{\alpha_1}\cdots\xi_d^{\alpha_d}\otimes \xi_{i_1}\wedge\cdots\wedge\xi_{i_p}&
\longmapsto &
{1\over {n}}\sum_{m=1}^{p}(-1)^{m+1}\xi_1^{\alpha_1}\cdots\xi_d^{\alpha_d}\xi_{i_m}\otimes
\xi_{i_1}\wedge\cdots\wedge\widehat{\xi}_{i_m}\wedge\cdots\wedge\xi_{i_p}.\cr
}
$$
We note that this definition also makes sense in the divided powers setting.
Indeed, in characteristic p we replace
$\xi_i^{\alpha_i}$ by $\xi_i^{[\alpha_i]}$ and the local description becomes:
$$
\matrix{
\c I^{n-p}/\c I^{n-p+1}\otimes_{\c O_X}\cOmega^p_X & \buildrel{s_p}\over\lra &
\c I^{n-p+1}/\c I^{n-p+2}\otimes_{\c O_X}\cOmega^{p-1}_X\hfill\cr
\xi_1^{[\alpha_1]}\cdots\xi_d^{[\alpha_d]}\otimes \xi_{i_1}\wedge\cdots\wedge\xi_{i_p}&
\longmapsto &
\sum_{m=1}^{p}(-1)^{m+1}\xi_1^{[\alpha_1]}\cdots\xi_{i_m}^{[\alpha_{i_m}+1]}\xi_d^{[\alpha_d]}\otimes
\xi_{i_1}\wedge\cdots\wedge\widehat{\xi}_{i_m}\wedge\cdots\wedge\xi_{i_p}.\cr
}
$$
Now let us compute the composition $D^{p-1}\opp s_p+ s_{p+1}\opp D^p$ on an element of the basis. We have
$$
\eqalign{
(D^{p-1}\opp s_p+& s_{p+1}\opp D^p)
(\xi_1^{\alpha_1}\cdots\xi_d^{\alpha_d}\otimes \xi_{i_1}\wedge\cdots\wedge\xi_{i_p})=\cr
=&D^{p-1}({1\over {n}}\sum_{m=1}^{p}(-1)^{m+1}\xi_1^{\alpha_1}\cdots\xi_d^{\alpha_d}\xi_{i_m}\otimes
\xi_{i_1}\wedge\cdots\wedge\widehat{\xi}_{i_m}\wedge\cdots\wedge\xi_{i_p})+\cr
&+s_{p+1}(\sum_{j=1}^{d}\alpha_j\xi_1^{\alpha_1}\cdots\xi_j^{\alpha_j-1}\cdots\xi_d^{\alpha_d}\otimes
\xi_j\wedge\xi_{i_1}\wedge\cdots\wedge\xi_{i_p})=\cr
=&
{1\over {n}}(\sum_{m=1}^{p}\sum_{j=1,\ne i_m}^{d}(-1)^{m+1}\alpha_j
\xi_1^{\alpha_1}\cdots\xi_j^{\alpha_j-1}\cdots\xi_d^{\alpha_d}\xi_{i_m}\otimes\xi_j\wedge
\xi_{i_1}\wedge\cdots\wedge\widehat{\xi}_{i_m}\wedge\cdots\wedge\xi_{i_p})+\cr
&+{1\over {n}}(\sum_{m=1}^{p}(\alpha_{i_m}+1)
\xi_1^{\alpha_1}\cdots\xi_d^{\alpha_d}\otimes \xi_{i_1}\wedge\cdots\wedge\xi_{i_p})+\cr
&+{1\over {n}}(\sum_{m=1}^{p}\sum_{j=1,\ne i_m}^{d}(-1)^{m}\alpha_j
\xi_1^{\alpha_1}\cdots\xi_j^{\alpha_j-1}\cdots\xi_d^{\alpha_d}\xi_{i_m}\otimes\xi_j\wedge
\xi_{i_1}\wedge\cdots\wedge\widehat{\xi}_{i_m}\wedge\cdots\wedge\xi_{i_p})+\cr
&+{1\over {n}}(\sum_{j=1\ne i_1,\dots , i_p}^{d}\alpha_{j}
\xi_1^{\alpha_1}\cdots\xi_d^{\alpha_d}\otimes \xi_{i_1}\wedge\cdots\wedge\xi_{i_p})=\cr
=& {1\over {n}}(p+\sum_{j=1}^d \alpha_j)
\xi_1^{\alpha_1}\cdots\xi_d^{\alpha_d}\otimes \xi_{i_1}\wedge\cdots\wedge\xi_{i_p}=\cr
=& \xi_1^{\alpha_1}\cdots\xi_d^{\alpha_d}\otimes \xi_{i_1}\wedge\cdots\wedge\xi_{i_p}\cr
}
$$
as desired.
}
\smallskip

The complex
$\{\c P^m_X\}_{\b Z}\otimes_{\c O_X}\cOmega^\point_X\cong
\ov{\DR}_X(\{\c P^m_X\}_{\b Z})$
where we take $\{\c P^m_X\}_{\b Z}=p_1^\ast(\c O_X)$. 
Thus we obtain a ``Poincar\'e Lemma":
$\c O_X\cong \ov{\DR}_X(\{\c P^m_X\}_{\b Z})$.
This result is the ``dual" of the quasi-isomorphism
$\cOmega^\point_X\otimes_{\c O_X}\c D_X\cong \omega_X[d]$ [B; VI.3.5].
\endgraf
The same result holds true on taking
$\{\c P^m_X\}_{\b Z}=p_0^\ast(\c O_X)$.
In this case we obtain the following corollary.

\smallskip

\coro{NN}
{
The complex
$$
\xymatrix{
\c O_X\ar[r]^{d_1}& \{\c P^m_X\}_{\b Z} \ar[r]&\cOmega^1_X\otimes_{\c O_X}\{\c P^m_X\}_{\b Z}\ar[r]&
\cdots \ar[r] &\cOmega^d_X\otimes_{\c O_X}\{\c P^m_X\}_{\b Z}\cr
}\leqno{\bf (\NNN)}\rif{tPoincris}
$$
is exact.  Hence,
$$
\c O_X \lra \ov{\DR}(p_0^\ast(\c O_X))
$$
is a quasi-isomorphism in $C^b(\nu(\c P^\cdot_X))$.
}
\smallskip

\coro{NN}
{\rif{cor}
Let $\c M\in \nu(\c P^\cdot_X)$.
The complexes
$$\xymatrix{
\c M\ar[r]^(0.3){d_1}& \{\c P^m_X\}_{\b Z}\otimes_{\c O_X} \c M\ar[r]&
\cOmega^1_X\otimes_{\c O_X}\{\c P^m_X\}_{\b Z}\otimes_{\c O_X} \c M\ar[r]&
\cdots \ar[r] &\cOmega^d_X\otimes_{\c O_X}\{\c P^m_X\}_{\b Z}\otimes_{\c O_X} \c M\cr
}\leqno{\bf (\NNN)}
$$
and
$$\xymatrix{
\c M\ar[r]^(0.3){d_0}& \c M\otimes_{\c O_X} \{\c P^m_X\}_{\b Z}\ar[r]&
\c M\otimes_{\c O_X}\{\c P^m_X\}_{\b Z}\otimes_{\c O_X}\cOmega^1_X \ar[r]&
\cdots \ar[r] &\c M\otimes_{\c O_X}\{\c P^m_X\}_{\b Z}\otimes_{\c O_X}\cOmega^d_X\cr
}\leqno{\bf (\NNN)}\rif{2}
$$
are exact.
}
\proof{
Let us consider the sequence (\cite{2}).
Its analogue for $\c M=\c O_X$
is the sequence \cite{Poincris}, which is exact and so locally homotopic to zero.
Any additive functor respects homotopies, 
so \cite{Poincris} tensorized over $\c O_X$ with $\c M$
gives \cite{2} which is exact.

}
\smallskip

\section{Equivalences of Categories}

\rema{NN}
{
Let $\c F^\point$ be an object in $C_1(\nu(\c O_X),\Diff_X)$.
By definition the differential
$d^i_{\c F}:\c F^i\lra \c F^{i+1}$ is a morphism of Pro-object
represented by  differential
operators of order at most one.
It defines in a unique way a morphism
$\ov d^i_{\c F}: \c P^1_X\otimes_{\c O_X}\c F^i\lra \c F^{i+1}$.
This morphism extends in a unique way to a morphism of stratified Pro-modules
$\Q^0_X( d^i_{\c F}): 
\{\c P_X^m\}_{\Bbb Z}\otimes_{\c O_X}\c F^{i}\lra\{\c P_X^m\}_{\Bbb Z}\otimes_{\c O_X} \c F^{i+1}$.
Using the $\c O_X$-base of $\c P_{X,1}$ given in local coordinates by 
$\b I, dx_1,\cdots ,dx_d$ we have for any section $s$ of $\c F^i$
$$
\eqalign{
d^i_{\c F}(s)=  & \ov d^i_{\c F}(\b I\otimes s)
\cr
d^i_{x_j}(s):= &\ov d^i_{\c F}(\xi_j\otimes s)\cr}
\leqno(\NNN)\rif{eq1}
$$
(where the second is taken as definition of $d^i_{x_j}$).
The maps $d^i_{x_j}:\c F^i\lra \c F^{i+1}$ are maps
of $\c O_X$-modules depending on the choice of the coordinates.
}
\smallskip

\defi{NN}
{
Let $\c F^\point\in C^b(\nu(\c O_X),\Diff_X)$.
Then $\c F^\point$ is a 
$\Omega_X^\point$-module in Pro-objects.
Let by definition
$\sigma ^{i,j}_{\c F}:\Omega^{j}_X\otimes_{\c O_X}
\c F^{j-i}\lra \c F^{i}$ 
be the $\Omega_X^\point$-structural maps.
They are:
$$
\sigma^{i,0}_{\c F}=id_{\c F^i}:  \c F^i  \lra \c F^i
$$
$$
\xymatrix{
\c P^1_X\otimes_{\c O_X}\c F^{i-1} 
\ar[r]^(.7){\ov{d}^{i-1}_{\c F^\point}} &\c F^i \\
\Omega_X^1\otimes_{\c O_X}\c F^{i-1}
\ar[u]\ar[ur]_{\sigma^{i,1}_{\c F}}\\
}
$$
and in general
$$
\xymatrix{
\c P^1_X\otimes_{\c O_X}\cdots \otimes_{\c O_X}\c P^1_X\otimes_{\c O_X}\c F^{i-j}
\ar[r]  &\cdots \ar[r]^(0.3){id_{\c P^1_X}\otimes \ov{d}^{i-2}_{\c F}}&
\c P^1_X\otimes_{\c O_X}\c F^{i-1}\ar[r]^(0.7){\ov{d}^{i-1}_{\c F}} &
\c F^i \\
\Omega^j_X\otimes_{\c O_X}\c F^{i-j}
\ar[u]\ar[urrr]_{\sigma^{i,j}_{\c F}} 	\\
}
$$
where the last vertical map 
$$
\Omega^j_X\otimes_{\c O_X}\c F^{i-j} \lra \overbrace{\c P^1_X\otimes_{\c O_X}\c P^1_X\otimes_{\c O_X}\cdots \otimes_{\c O_X}
\c P^1_X}^{\hbox{j times}}\otimes_{\c O_X}\c F^{i-j} 
$$
is induced by the shuffle on a basis.
}
\smallskip

Now we prove a technical lemma which will be used in the proof of our main theorem.

\lemm{NN}
{\rif{condiff}
Given $\c F^\point\in C_1(\nu(\c O_X),\Diff_X)$; the morphisms $d^i_{\c F}$, 
$d^i_{x_j}$ of (\cite{eq1}) 
for $i\in \Bbb Z$ and $j\in \{0,...,d\}$
satisfy the following conditions:
	\item{i)} $d^{i+1}_{\c F}\opp d^{i}_{\c F}=0$
	\item{ii)} $d^{i+1}_{x_j}\opp d^i_{\c F}+
	            d^{i+1}_{\c F}\opp d^i_{x_j}=0$
	\item{iii)}$d^{i+1}_{x_j}\opp d^i_{x_k}+
		    d^{i+1}_{x_k}\opp d^i_{x_j}=0$	
	\item{iv)} $d^{i+1}_{x_j}\opp d^i_{x_j}=0$.
}
\proof{
The first condition is given by the hypothesis $\c F^\point\in
C_1(\nu(\c O_X),\Diff_X)$, so the composition
$d^{i+1}_{\c F}\opp d^{i}_{\c F}=0$.
In order to prove ii) to iv) we remark that:
$$
\eqalign{
d^{i+1}_{\c F}\opp d^i_{x_j}(s)+d^{i+1}_{x_j}\opp d^i_{\c F}(s)= &
\ov{d^{i+1}_{\c F} d^i_{\c F}}(\xi_j\otimes s)=0 \cr
d^{i+1}_{x_j}\opp d^i_{x_k}(s)+d^{i+1}_{x_k}\opp d^i_{x_j}(s)= &
\ov{d^{i+1}_{\c F} d^i_{\c F}}(\xi_j \xi_k\otimes s)=0\cr
d^{i+1}_{x_j}\opp d^i_{x_j}(s)+d^{i+1}_{x_j}\opp d^i_{\c F}(s)= &
\ov{d^{i+1}_{\c F}d^i_{\c F}}(\xi^2_j\otimes s)=0 \cr
}$$
where $\xi_j=1\otimes x_j-x_j\otimes 1\in \c P^1_X$.
We recall that $\ov{d^{i+1}_{\c F}d^i_{\c F}}$ is obtained
by the composition
$$\xymatrix{
\c P^1_X\otimes_{\c O_X}\c P^1_X\otimes_{\c O_X}\c F^i 
\ar[r]^(0.6){id_{\c P^1_X}\otimes \ov{d^i_{\c F}}} & 
\c P_X^1\otimes_{\c O_X}\c F^{i+1} 
\ar[r]^(0.6){\ov{d^{i+1}_{\c F}}} &
\c F^{i+2} \cr
\c P^2_X\otimes_{\c O_X}\c F^i
\ar[u]^{\delta^{2,1}\otimes id_{\c F^i}} \ar[urr]_{\ov{d^{i+1}_{\c F}d^i_{\c F}}} & &\cr
}$$
and $\delta^{2,1}$ was defined in  \cite{localbase}.
}\smallskip

\defi{NN}{
Let $dx_{1}, \dots, dx_n$ be a local basis for $\Omega^1_X$
(where $n$ is the dimension of $X$).
We define the maps
$$\matrix{
\eta^{i,j}_{\c F}:& \Omega^j_X\otimes_{\c O_X}\c F^{i-j} & \lra &  \c F^i\hfill \cr
& f dx_{i_1}\wedge \cdots \wedge dx_{i_j}\otimes s & \longmapsto &
f dx_{i_1}^{i-1}\opp \cdots \opp dx_{i_j}^{i-j}(s) \cr}
$$
where $f$ is a section of $\c O_X$ and $s$ is a section of $\c F^{i-j}$.
This definition does not depend on local coordinates and moreover 
$\sigma^{i,j}_{\c F}=j!\eta^{i,j}_{\c F} $.
}
\smallskip

\defi{NN}
{\rif{defiphi}
Let $q=\{q_m\}_{\Bbb Z}:\{\c P_X^m\}_{\Bbb Z}\lra \c O_X$ be the usual projection which is linear for
both the $\c O_X$-module structures of $\{\c P^m_{X}\}_{\Bbb Z}$.
Given $\c F^{\point}\in C_1(\nu(\c O_X),\Diff_X)$
we define the morphisms
$$
\Phi^i_{\c F}: \bigoplus_{j=0}^d 
\Omega^j_X\otimes_{\c O_X}\{\c P_X^m\}_{\Bbb Z}\otimes_{\c O_X}\c F^{i-j}
\lra 
\c F^i
$$
for each $i\in \Bbb Z$
in the following way:
we consider the composition
$$
\xymatrix{
\Omega^j_X\otimes_{\c O_X}\{\c P_X^m\}_{\Bbb Z}\otimes_{\c O_X}\c F^{i-j}
\ar[r]^(.6){id_{\sss{\Omega^j_X}}\otimes q\otimes id_{\sss{\c F^{i-j}_X}}}
\ar[dr]_(.6){\Phi^{i,j}_{\c F}}	&
\Omega^j_X\otimes_{\c O_X}\c F^{i-j}
\ar[d]^{\eta^{i,j}_{\c F}}	\\								&
\c F^{i}	\\
}
$$
and by definition
$\Phi^i_{\c F}:=\sum_{j=0}^d \Phi^{i,j}_{\c F}$.
}
\smallskip

%

\theo{NN}
{\rif{bigtheo}
We have two morphisms of functors
$$
\Psi: id_{C^b(\c P^\cdot_X)}
 \lra {\Q}^0_X \opp \ov{\DR}_X=
  {\Q}^0_{1,X} \opp \ov{\DR}_{1,X}
 $$
(functors between $C^b(\c P^\cdot_X)$
and itself)
$$
\Phi:\ov{\DR}_{1,X}\opp {\Q}^0_X\lra id_{C^b_1(\nu(\c O_X),\Diff_X)}
$$
(functors between $C^b_1(\nu(\c O_X),\Diff_X)$
and itself).
They induce  quasi-isomorphisms
of complexes.
So the functors $\ov{\DR}_{1,X}$ and $\ov{\DR}_X$ localize with respect to
${\Q}^0_X$-quasi-isomorphisms
inducing the functor
$$
\ov{\DR}_{1,X}:
D^b(\c P^\cdot_X)\lra D^b_1(\nu(\c O_X),\Diff_X).
$$
$$
\ov{\DR}_X:
D^b(\c P^\cdot_X)\lra D^b(\nu(\c O_X),\Diff_X).
$$
Moreover  $\ov{\DR}_X$ (resp. $\ov{\DR}_{1,X}$) is an equivalence of categories
whose quasi-inverse is the functor
${\Q}^0_X$ (resp. ${\Q}^0_{1,X}$).
}
\proof{
The morphism
$d^0:\c O_X\ra \c P^m_X$ induces a morphism of bicomplexes (see (\cite{2}))
$$
\Psi_{\c M^\point}: \c M^\point \lra
\c M^\point\otimes_{\c O_X}\{\c P^m_X\}_{\b Z}\otimes_{\c O_X}\cOmega^\point_X.
$$
Then we obtain a morphism of complexes
$$\eqalign{
\c M^\point\lra &
\c M^\point\otimes^\point_{\c O_X}\{\c P^m_X\}_{\b Z}\otimes_{\c O_X}\cOmega^\point_X\cong\cr
\cong&  \{\c P^m_X\}_{\b Z}\otimes_{\c O_X}\c M^\point\otimes^\point_{\c O_X}\cOmega^\point_X\cong\cr
\cong&\Q^0_X(\ov{\DR}_X(\c M^\point))\cong\cr
\cong&\Q^0_{1,X}(\ov{\DR}_{1,X}(\c M^\point)).\cr
}
$$
The isomorphism between the first and the second complex  is induced
by the stratification on $\c M^\point$ which is the isomorphism 
$ \{\c P^m_X\}_{\b Z}\otimes_{\c O_X}\c M^\point\cong
\c M^\point\otimes_{\c O_X}\{\c P^m_X\}_{\b Z}$.
By Corollary \cite{cor} we obtain
that it is a quasi-isomorphism.
So the functors $\ov{\DR}_{1,X}$ and $\ov{\DR}_X$ send qis in
$D^b(\c P^\point_X)$ into $\Q^0_X$-qis which permits
us to localize them.
\endgraf
We have to prove that
the diagram
$$
\xymatrix{
\Omega^j_X\otimes_{\c O_X}\{\c P_X^m\}_{\Bbb Z}\otimes_{\c O_X}\c F^{i-j}
\ar[r]^{d^i_{\ov{\DR}\Q^0}} \ar[d]_{\Phi^i_{\c F}} &
\Omega^j_X\otimes_{\c O_X}\{\c P_X^m\}_{\Bbb Z}\otimes_{\c O_X}\c F^{i+1-j}
\ar[d]^{\Phi^{i+1}_{\c F}} \\
\c F^i \ar[r]_{d^i_{\c F}} & \c F^{i+1} \\
}
$$
is commutative.
\endgraf
The complex
$
{\ov{\DR}_X}{\Q_X}^{0}(\c F^{\point})=
(\c G^{\point\point})_{tot}
$
where
$$
\c G^{p,q}=
\Omega^p_X\otimes_{\c O_X}\{\c P_X^m\}_{\Bbb Z}\otimes_{\c O_X}\c F^{q}
$$
and
$$
d^{\prime p,q}_{\c G}:
\Omega^{p}_X\otimes_{\c O_X}\{\c P_X^m\}_{\Bbb Z}\otimes_{\c O_X}\c F^{q}
\lra  
\Omega^{p+1}_X\otimes_{\c O_X}\{\c P_X^m\}_{\Bbb Z}\otimes_{\c O_X}\c F^{q}$$
is
$$
d^{\prime p,q}_{\c G}=					 
d^p_{\ov{\DR}_X(p^\ast_0(\c O_X))}\otimes id_{\c F^q}
\leqno{\bf (\NNN)}\rif{d1}
$$
while
$$
d^{\prime\prime p,q}_{\c G}:
\Omega^{p}_X\otimes_{\c O_X}\{\c P_X^m\}_{\Bbb Z}\otimes_{\c O_X}\c F^{q}
\lra  
\Omega^{p}_X\otimes_{\c O_X}\{\c P_X^m\}_{\Bbb Z}\otimes_{\c O_X}\c F^{q+1}
$$ 
is 
$$
d^{\prime\prime p,q}_{\c G}=	
id_{\Omega_X^p}\otimes \Q^0_X(d^q_{\c F}) ; 	
\leqno{\bf (\NNN)}\rif{d2}
$$
where 
$\ov{\DR}_X(p^\ast_0(\c O_X))$ was considered in  \cite{tPoincris}.
\endgraf

Given $I^{\point,\point}$ 
a bounded bicomplex with commuting differentials 
$d^{\prime p,q}_{I}:I^{p,q}\ra I^{p+1,q}$ and 
$d^{\prime\prime p,q}_{I}:I^{p,q}\ra I^{p,q+1}$,
the total complex associated to it is denoted by $I_{\tot}^\point$
with
$
I_{\tot}^{r}\bydef 
\bigoplus _{p+q=r}I^{p,q}
$
and
$
d_{I_{\tot}}(x)=d^{\prime p,q}_I(x)+(-1)^p d^{\prime\prime p,q}_I(x)
$
for any $x\in I^{p,q}$. \noindent
\par
The set 
$\Hom_{} (A_{\tot}^{\point},B^{\point})$
of morphisms of complexes 
between a total complex of a bicomplex
and a complex is the set families of maps
$\{\phi^p:A_{\tot}^p\ra B^p\}_p$
such that 
$d^p_{B}\circ \phi^p=\phi^{p+1}\circ d^p_{A_{\tot}}$.
Then $\Hom_{} (A_{\tot}^{\point},B^{\point})$ is
isomorphic to the set of families of maps 
 $\{\phi^{p,q}:A^{q,p-q}\ra B^{p}\}_{p,q}$ 
 satisfying the following conditions 
 $$ 
d^p_B\circ \phi^{p,q}=
\phi^{p+1,q+1}d_{A}^{\prime q,p-q}
+(-1)^{q}\phi^{p+1,q}\circ d_{A}^{\prime\prime q,p-q} 
\leqno{\bf (\NNN)}\rif{conmortot}$$ 
for any $p,q$.
So we have only to prove that
$$
d^p_{\c F}\circ \Phi^{p,q}_{\c F}=
\Phi^{p+1,q+1}_{\c F}d_{\c G}^{\prime q,p-q}
+(-1)^{q}\Phi^{p+1,q}_{\c F}\circ d_{\c G}^{\prime\prime q,p-q} 
 \leqno{\bf (\NNN)}\rif{conmortot}
$$
is true.\endgraf
It is enough to check these relations locally, choosing local coordinates
$x_1,\dots, x_n$.
The sections $f dx_{i_1}\wedge \cdots \wedge dx_{i_q}\otimes \b I\otimes s$
generate $\Omega^{q}_X\otimes_{\c O_X}\{\c P_X^m\}_{\Bbb Z}\otimes_{\c O_X}\c F^{p-q}$
where $s$ is a section of $\c F^{p-q}$. Then
$$\eqalign{
d^p_{\c F}\circ \Phi^{p,q}_{\c F}
(f dx_{i_1}\wedge \cdots \wedge dx_{i_q}\otimes \b I\otimes s)=&
d^p_{\c F}(f d^{p-1}_{x_{i_1}}\cdots d^{p-q}_{x_{i_q}}(s))=\cr
= &f d^p_{\c F}d^{p-1}_{x_{i_1}}\cdots d^{p-q}_{x_{i_q}}(s)+
\sum_{i=1}^n {\partial f\over\partial x_i} d^p_{x_i}
d^{p-1}_{x_{i_1}}\cdots d^{p-q}_{x_{i_q}}(s) \cr
}$$
while
$$\eqalign{
\Phi_{\c F}^{p+1,q+1} d^{\prime q,p-q}_{\c G}
(f dx_{i_1}\wedge \cdots \wedge dx_{i_q}\otimes \b I\otimes s)=&
\Phi_{\c F}^{p+1,q+1}(\sum_{i=1}^n{\partial f\over\partial x_i}dx_i\wedge 
dx_{i_1}\wedge\cdots\wedge dx_{i_q}\otimes \b I\otimes s)=\cr
=& \sum_{i=1}^n{\partial f\over\partial x_i}
d^p_{x_i}d^{p-1}_{x_{i_1}}\cdots d^{p-q}_{x_{i_q}}(s)
}
$$
For the last term we have
$$
\eqalign
{
\Phi^{p+1,q}_{\c F}d^{\prime\prime q,p-q}_{\c G}
(f dx_{i_1}\wedge \cdots \wedge dx_{i_q}\otimes \b I\otimes s)=&
\Phi^{p+1,q}_{\c F}(f dx_{i_1}\wedge \cdots 
\wedge dx_{i_q}\otimes \b I\otimes d^{p-q}_{\c F}(s))=\cr
=& f d^p_{x_{i_1}}\cdots d^{p-q+1}_{x_{i_q}}d^{p-q}_{\c F}(s)\cr
}
$$
Thus, using Lemma \cite{condiff},
we prove our assertion.
Moreover the composition 
$$\xymatrix{
\c F^\point\ar[r]^(0.3){d^1}\ar[dr]_{id_{\c F^\point}}& 
\cOmega^\point_X\otimes^\point_{\c O_X}\{\c P^m_X\}_{\b Z}\otimes_{\c O_X}
\c F^\point\ar[d]_{\Phi^\point_{\c F}}\cr
&\c F^\point}
$$
is the identity so $\Phi_{\c F^\point}$ is a $\Q^0_X$-quasi-isomorphism.
In particular the functor $\ov{\DR}_X$ localizes to
$\ov{\DR}_X:D^b(\c P^\cdot_X)\lra D^b(\nu(\c O_X),\Diff_X)$
and it is an equivalence of categories with quasi-inverse the localized
$\Q^0_X$ functor.
}
\smallskip

\coro{NN}
{
The functor
$$\lambda: D^b_1(\nu(\c O_X),\Diff_X)\lra D^b(\nu(\c O_X),\Diff_X)
$$ is an equivalence of categories whose quasi-inverse is the functor
$\ov{\DR}_X\opp\Q^0_X$.
}
\smallskip

%


%

\section{Crystals in Pro-Modules}
We refer to Grothendieck expos\'e [G] 
for the definition of crystalline site $Cris(X/\b C)$ in characteristic zero (and also
to Berthelot's thesis [Be]).
We denote by $\c O_{X/\b C}$ the sheaf on $Cris(X/\b C)$
such that for any object (nilpotent closed immersion
$U\hookrightarrow T$
with $U\subset X$ open subset) its value is 
$\c O_{X/\b C}(U\hookrightarrow T)\bydef \c O_T$.
It is a ringed sheaf on  $Cris(X/\b C)$.

\defi{NN}
{
A crystal in Pro-modules $\{\c F_i\}_{i\in I}$
is a sheaf on Pro-$\c O_{X/\b C}$-modules in the crystalline site
$Cris(X/\b C)$ such that for any morphism 
$p:\{U'\hookrightarrow T'\}\ra \{U\hookrightarrow T\}$
given by the diagram
$$
\xymatrix{
U'\ar[r]\ar[d] & T' \ar[d]_p \cr
U\ar[r]\ar[d] & T\cr
X
 \cr
}
$$
we have
$\{p^{\ast}(\c F_i(U\hookrightarrow T))\}_{i\in I}\cong
\{\c F_i(U'\hookrightarrow T')\}_{i\in I}$.
}
\smallskip

\theo{NN}
{
The category
$\nu(\c P^\cdot_X)$ is equivalent to the category of
crystals in Pro-coherent $\c O_{X/\b C}$-modules.
}
\proof{
The proof is equivalent to the classical proof of the equivalence
between stratified $\c O_{X/\b C}$-modules and crystals
(see [Be]).
If $\{\c F_i\}_{i\in I}$ is a crystal in Pro-coherent $\c O_{X/\b C}$-modules
 then
$\{\c F_i(X\buildrel{id_X}\over\ra X)\}_{i\in I}\in \nu(\c O_X)$.
Taking the diagram defined by the diagonal thickenings 
$X_n\lra X\times X$ with $n\in \b N$
$$
\xymatrix{
X\ar[r] \ar[d]_{id}& X_n\ar[d]_{p_0}\ar@<1ex>[d]^{p_1}\cr
X\ar[r]_{id}&  X \cr
}
$$
we obtain the stratification
$\{p_0^{\ast}(\c F_i(X\hookrightarrow X_n))\}_{i\in I}\cong
\{p_1^{\ast}(\c F_i(X\hookrightarrow X_n))\}_{i\in I}$
in the Pro-object $\{\c F_i(X\buildrel{id_X}\over\ra X)\}_{i\in I}$.
\endgraf\noindent
Conversely if $\{\c F_i\}_I\in \nu(\{\c P^m_X\}_{\b Z})$
we define a sheaf on $Cris(X/\b C)$ in the following way.
For any object $U\hookrightarrow T\in Cris(X/\b C)$
there exists locally a section $h:T\ra X$ ($X$ is smooth).
 We define
$CR(\{\c F_i\}_I)(U\hookrightarrow T)\bydef 
\{h^{\ast}(\c F_i)\}_I$.
These local definitions patch together to define a sheaf
in Pro-coherent $\c O_{X/\b C}$-Modules which is a crystal.  
}
\smallskip


\rema{NN}
{
Deligne proved in a conference at IHES (1970 unpublished)
that the category of ``regular" crystals in Pro-modules on $X$ is equivalent 
to the category of "algebraic" constructible sheaves in the analytic space
$X^{an}$.
Hence, Theorem \cite{bigtheo} might be interpreted 
 as an algebraic Pro-version of a Riemann-Hilbert
equivalence.
}
\smallskip

\chiudiriferimenti

\refer
\biblio{AB}
                    {Andr\'e Y. and Baldassarri  F.
                     {De Rham Cohomology of Differential Modules on Algebraic
			Varieties.}
                     {\it Progress in Mathematics},
                     {\bf 189} Birkh\"auser Verlag, Basel, 2001.}
                     
\biblio{AM}
       {Artin M., Mazur B.,
        {\'Etale Homotopy. }
        {Lectures Notes in Mathematics, Vol.100.,  }
        {\it  Springer-Verlag, Berlin-New York,} 1969 }

\biblio{B}
                    {Borel et al.
                     {\it Algebraic D-modules.}
                     Perspectives in Mathematics, Vol. 2
			J. Coates and S. Helgasan editors.}

\biblio{Be}
                    {Berthelot  P.
                     {Cohomologie cristalline des sch\'emas de
			caract\'eristique $p>0$.}
                     {Lecture Notes in Mathematics, Vol. 407.}
                     {\it Springer-Verlag, Berlin-New York,} 1974.}


\biblio{BCF}
                    {Baldassarri  F., Cailotto M. and Fiorot L.
                     {Poincar\'e Duality for Algebraic De Rham
                     Cohomology.},
                     {\it Manuscripta Math.}
			{\bf 114},1 (2004)  61-116.}

\biblio{BeO}
       {Berthelot P and Ogus A. 
        {Notes on crystalline cohomology.}
        {\it Princeton University Press, Princeton, N.J.; 
           University of Tokyo Press, Tokyo,} 
        1978.}


\biblio{DB}
                    {Du Bois Ph.
                     {Complexe de de Rham filtr\'e d'une vari\'et\'e
			singuli\`ere.}
                     {\it Bull. Soc. Math. France},
                     {\bf 109} (1981), 41--81.}


\biblio{EGA}
                    {Grothendieck A. and Dieudonn\'e J.
                     El\'ements de g\'eom\'etrie alg\'ebrique.
                     {\it Inst. Hautes \'Etudes Sci. Publ. Math.},
                     {\bf 4} (1960), {\bf 8} (1961), {\bf 11} (1961), {\bf
			17} (1963),
                     {\bf 20} (1964), {\bf 24} (1965), {\bf 28} (1966), {\bf
			32} (1967).}

\biblio{F.1} {Fiorot L.
                     {On derived category of differential complexes.},
                     Journal of Algebra 312 (2007), 362--376  }

\biblio{F.2} {Fiorot L.
                     {A simple introduction to Ind and Pro-categories.}      
                     }

\biblio{G}
                    {Grothendieck A. 
                     {Crystals and De Rham Cohomology of Schemes}
                     {\it Dix expos\'es sur la cohomologie des
                     sch\'emas}, 
			{\it North-Holland, Amsterdam,} 1968}

\biblio{H.1}
   {Hartshorne R.
    {Cohomology with compact supports for coherent
	 sheaves on an algebraic variety.}
    {\it Math. Ann.},
    {\bf 195} (1972), 199--207.}

\biblio{H.2}
   {Hartshorne R.
    {On the De Rham cohomology of algebraic varieties.}
    {\it Inst. Hautes \'Etudes Sci. Publ. Math.},
    {\bf 45} (1975), 5--99.}

\biblio{HL}
                    {Herrera M. and Lieberman D.
                     {Duality and the de Rham cohomology of infinitesimal
			neighborhoods.}
                     {\it Invent. Math.},
                     {\bf 13} (1971), 97--124.}

\biblio{H.RD}
                    {Hartshorne R.
                     {\it Residues and duality.}
                     {Lecture Notes in Mathematics, Vol. 20}}
                     
 \biblio{J}
                    {Jannsen U.
                     {Continuous \'etale cohomology.}
                     {\it Math. Ann.}, {\bf 280} (1988), no. 2, 207--245.}

\biblio{L}
                    {Laumon G.
                     {Sur la cat\'egorie d\'eriv\'ee des $\c D$-modules filtr\'es.}
                     {\it Algebraic Geometry}, proceedings, Tokyo/Kyoto 1982
                     Lecture Notes in Mathematics, 1016, Springer-Verlag.}

\biblio{M}
                    {Mebkhout Z.
                     {\it Le formalisme des six op\'erations de Grothendieck
			pour les $\c D\sb X$-modules coh\'erents.}
                     {Travaux en Cours, 35},
                     {\it Hermann, Paris,} 1989.}

\biblio{S.1}
                    {Saito M.
                     {Modules de Hodge polarisables.}
                     {\it Publ. RIMS, Kyoto Univ.},
                     {\bf 24} (1988), 849--995.}

\biblio{S.2}
                    {Saito M.
                     {Induced $\c D$-modules and differential complexes.}
                     {\it Bull. Soc. Math. France},
                     {\bf 117} (1989), 361--387.}

\bigskip

\end

In fact using the commutativity of (\cite{comdiadiop})
we obtain a morphism of bicomplexes
$$
\c M^\point\buildrel{d^1\otimes id_{\c M^\point}}\over\lra
\{\c P^m_X\}_{\b Z}\otimes_{\c O_X}(\c M^{\point}\otimes_{\c O_X}\cOmega^\point_X)
$$
so a map
$$
\c M^\point\lra {\Q}^0_X[-d]\opp \ov{\DR}_X(\c M^\point).
$$

\defi{NN} 
{ 
Let  
$\cOmega^{\point}_{X}$  
be the complex of relative differential forms, 
$\cTheta^{\point}_{X}$  
be the relative tangent complex  
(so $\cOmega^i_{X}=\cHom_{\c O_X}(\cTheta^{-i}_{X},\c O_X)$)  
and 
$\omega_{X}=\cOmega^{d_{X}}_{X}$ where  
$d_{X}=\dim X=\dim X-\dim S$ is the relative dimension. 
} 
\bigskip 

 
\defi{NN} 
{\rif{d1} [EGA IV, 4.16.7]
Let  
$ 
\Delta : X \ra X \times_{S} X 
$ 
be the diagonal embedding.  
The Pro-object 
$\{ \c P_{X}^{m}\}_{m\in \Bbb Z}$  
is defined as  
$
\c P_{X}^{m}:=
\Delta^{-1}(\c O_{X\times_{S} X}/\c I^{m+1})
$  
if $m\geq 0$
and zero otherwise
where $\c I$ is the ideal of the locally closed immersion 
$\Delta$. 
The transition maps are called
$
q_{n,m}: \c P_{X}^{n}\lra \c P_{X}^{m}
$
where $m\leq n$ and they are epimorphisms.  
For each $m\in\Bbb N$ we call 
$
q_m:=q_{m,0}:\c P_{X}^{m}\ra \c O_X
$.  
} 
\bigskip 
 
\rema{NN} 
{ 
Every $\c P_{X}^{m}$ has two canonical structures of  
$\c O_X$-module induced by the projections   
$
\pi_{i}: X\times_{S}X\ra X
$ 
($i=1,2$).  
\endgraf 
By notation every time we form a tensor product  
with $\c P_{X}^m$, 
we use the structure nearest to the tensor product:  
so for each $\c O_X$-module $\c F$,  
$\c P_{X}^{m}\otimes _{\c O_X}\c F$  
is made using the $\pi_2$ structure;  
on the other side $\c F\otimes_{\c O_X} \c P_{X}^{m}$  
uses the $\pi_1$ structure. 
\endgraf 
We remark that $p_m$ is $\c O_X$-linear for both these 
structures.  
} 
\bigskip

\defi{NN} 
{ 
We denote by 
$\Delta_\alpha:X\lra X \times_S \dots \times_S X$
the diagonal embedding ($\alpha+1$ factors).  
Let 
$\c P_{X}^\alpha(m):=\Delta^{-1}(\c O_{X\times_S \dots \times_S X}/
\c I^{m+1}_\alpha)$
(where $\c I_\alpha$ is the ideal of the locally closed immersion $\Delta_\alpha$).
\endgraf
There are $\alpha +1$ projections $\pi_i:X \times_S \dots \times_S X\lra X$
for $i=\{1,2,\dots,\alpha +1\}$ 
which induce $\alpha +1$ different $\c O_X$-module structures on $\c P_{X}^\alpha(m)$,
one for each projection.
We denote by the right structure the one induced by $\pi_{\alpha +1}$. 
\endgraf 
For any object  $\c F^\point$ in $C(\Coh(\c O_X),\Diff_{X})$ 
Grothendieck constructs 
functorially a bicomplex (using the right structure on
$\c P_{X}^\alpha(m)$) 
$$
{\Q}^{\point}_{X}(\c F^\point) 
=\{\c P_{X}^{\point}(m) \otimes_{\c O_X}\c F^\point\}_\b Z
$$ 
in the category 
$\Pro^{AR}(\Coh(\c O_X))$ 
such that the 
bicomplex 
$\limpr_m {\Q}^{\point}_{X}({\cal F}^{\point})$ 
is a resolution of 
$\c F^\point$. 
\endgraf
So we have a notion of ${\Q}^0_{X}$ functor. 
It is called the Grothendieck formalization. 
} 
\bs
 
 \defi{NNN} 
{ {\bf  DeRham functor for Stratified Modules.}\endgraf
Let consider the functor which sends an induced object in 
$\nu(\{\c P^m_{X}\}_{\Bbb Z})$
into his horizontal sections:
$$
\matrix{
\ov{\DR}_{X}: & \nu_i(\{\c P^m_{X}\}_{\Bbb Z})&\lra &
 M(\nu(\c O_X)),\Diff_{X})
\hfill\cr
& \{\c P^m_{X}\}_{\Bbb Z}\otimes_{\c O_X}\{\c E_i\}_I &
\longmapsto &\{\c E_i\}_I =
\cHom_{\{\c P^m_{X}\}_{\Bbb Z}}(\c O_X,
\{\c P^m_{X}\}_{\Bbb Z}\otimes_{\c O_X}\{\c E_i\}_I).\cr
}$$
When we restrict this functor to perfect complexes, which are induced by locally free
$\c O_X$-modules, the $\ov{\DR}_{X}$ takes image into the category
$M(\Coh(\c O_X),\Diff_{X})$.
So this functor induces a functor 
$$
\ov{\DR}_{X}:D^b_p(\{\c P^m_{X}\}_{\Bbb Z})\lra D^b_{S,p}(X).
$$
}
\bs

Let consider $\c P_I\in \nu(\c O_X)$ such that any $\c P_i$ is 
$\c O_X$-locally free of finite rank.
The functor
$\cHom_{\c O_X}(\c P_i,\_)$ is exact by hypothesis and
filtering inductive limits are exacts.
Any object $\c F_J\in \nu(\c O_X)$ can be represented as a Pro-object
whose  transition morphisms are surjective maps.
In this case the $\Pro(\QCoh(\c O_X))$-module
$\limind_I\cHom_{\c O_X}(\c P_i,\c F_J)$
is Mittag-Leffer.
Then the functor
$$
\cHom_{\nu(\c O_X)}(\c P_I,\_)\bydef
\limpro\limind_I\cHom_{\c O_X}(\c P_i,\_)
$$ 
is exact, which proves that
$\c P_I$ is locally projective.\par
Let us we prove that 
$\nu(\c O_X)$ has enough
locally projective.
Let $\c F_I\in \nu(\c O_X)$. 
For any $i\in I$
there exist a locally projective $\c O_X$-module
$\c P_i$ and an epimorphism
$\c P_i \buildrel{p_i}\over\lra \c F_i$.
Let us fix an index $i_0\in I$.
For any
$i\geq i_0$ let consider the diagram
$$
\xymatrix{
\c P_{ i_0}\ar[r]^{p_{i_0}}& \c F_{ i_0} \cr
&\c F_i\ar[u]\cr
&\c P_i \ar[u]_{p_i} .    \cr
}
$$
We can find a morphism
$\c P_i\lra \c P_{ i_0}$ which makes the diagram commutative, so we obtain a locally
projective Pro-coherent $\c P_{i\geq { i_0}}$.
\endgraf
Let us define a surjective map 
$$\psi\in \Hom_{\nu(\c O_X)}(\c P_{i\geq i_0},\c F_I)=
\limpro_I\limind_{i\geq i_0}\Hom_{\c O_X}(\c P_j,\c F_i)$$
in the following way: for $i\geq i_0$ we get $p_i:\c P_i\twoheadrightarrow \c F_i$
for $i\leq i_0$ we get $\pi_{\c F,i, i_0}\circ p_{i_0}:\c P_{i_0}\twoheadrightarrow \c F_i$
where $\pi_{\c F,i, i_0}:\c F_{i_0}\twoheadrightarrow \c F_i$ is the 
structural morphism of the projective system $\c F_I$.
For any other $i\in I$ let $j\in I$ be the first index in $I$ such that
$i\leq j, i_0\leq j$, then we get $\pi_{\c F,j i}\circ p_{j}:\c P_j\twoheadrightarrow \c F_i$.
\endgraf